\input ./style/arxiv-general.cfg
\documentclass[aap,MSNbibl,seceqn,dvips]{arximspdf}
\makeatletter
   \@ifpackageloaded{graphicx}{}{\usepackage{graphicx}}
\makeatother

\doi{10.1214/15-AAP1106}
\volume{26}
\issue{2}
\pubyear{2016}
\firstpage{857}
\lastpage{891}
\docsubty{FLA}

\makeatletter
\def\mid{|}
\newcommand{\eqref}[1]{(\ref{#1})}
\newtheorem{thmm}{Theorem}[section]
\newtheorem{lemma}[thmm]{Lemma}
\newtheorem{corollary}[thmm]{Corollary}
\newtheorem{prop}[thmm]{Proposition}
\newtheorem{propp}{Proposition}[section]

\newproclaim{defn}[thmm]{Definition}
\newproclaim{example}[thmm]{Example}
\newproclaim{examples}[thmm]{Examples}
\newproclaim{rem}[thmm]{Remark}
\newproclaim{question}[thmm]{Question}

\newcommand{\N}{\mathbb{N}}
\newcommand{\E}{\mathbb{E}}
\renewcommand{\P}{\mathbb{P}}

\newcommand{\calC}{\mathcal{C}}
\newcommand{\calL}{\mathcal{L}}

\newcommand{\Bin}{\operatorname{Bin}}
\newcommand{\Hyp}{\operatorname{Hyp}}

\makeatother

\begin{document}
\begin{frontmatter}

\title{A new coalescent for seed-bank models}
\runtitle{Seed-bank coalescent}

\begin{aug}
\author[A]{\fnms{Jochen} \snm{Blath}\ead[label=e1]{blath@math.tu-berlin.de}\thanksref{T1}},
\author[A]{\fnms{Adri\'an} \snm{Gonz\'alez Casanova}\ead[label=e2]{adriangcs@hotmail.com}\thanksref{T2}},\\
\author[A]{\fnms{Noemi} \snm{Kurt}\ead[label=e3]{kurt@math.tu-berlin.de}\thanksref{T1}}
\and
\author[A]{\fnms{Maite} \snm{Wilke-Berenguer}\corref{}\ead[label=e4]{wilkeber@math.tu-berlin.de}\thanksref{T3}}
\runauthor{Blath, Gonz\'alez Casanova, Kurt and Wilke-Berenguer}
\thankstext{T1}{Supported by DFG SPP 1590 ``Probabilistic
structures in evolution.''}
\thankstext{T2}{Supported by DFG RTG 1845
``Stochastic Analysis and Applications in Biology, Finance and
Physics,'' the Berlin Mathematical School (BMS), and the Mexican
Council of Science (CONACyT) in collaboration with the German Academic
Exchange Service (DAAD).}
\thankstext{T3}{Supported by DFG RTG 1845 ``Stochastic Analysis and Applications
in Biology, Finance and Physics'' and the Berlin Mathematical School (BMS).}
\affiliation{Technische Universit\"at Berlin}
\address[A]{Institut f\"ur Mathematik, Sekr. MA 7-5\\
Technische Universit\"at Berlin\\
Stra\ss e des 17. Juni 136\\
10623 Berlin\\
Germany\\
\printead{e1}\\
\phantom{E-mail:\ }\printead*{e2}\\
\phantom{E-mail:\ }\printead*{e3}\\
\phantom{E-mail:\ }\printead*{e4}}
\end{aug}

%
\received{\smonth{11} \syear{2014}}
%
\revised{\smonth{2} \syear{2015}}

%
\begin{abstract}
We identify a new natural coalescent structure, which we call the \emph
{seed-bank coalescent},
that describes the gene genealogy of populations under the influence of
a strong
seed-bank effect, where ``dormant forms'' of individuals (such as seeds
or spores) may jump a significant number of generations before joining
the ``active'' population.
Mathematically, our seed-bank coalescent appears as scaling limit in a
Wright--Fisher model with geometric seed-bank age structure if the
average time of seed dormancy scales with the order of the total
population size $N$. This extends earlier results of Kaj, Krone and
Lascoux [\textit{J. Appl. Probab.} \textbf{38} (2011) 285--300]
who show that the genealogy of a Wright--Fisher model in
the presence of a ``weak'' seed-bank effect is given by a suitably
time-changed Kingman coalescent.
The qualitatively new feature of the seed-bank coalescent is that
ancestral lineages are independently blocked at a certain rate
from taking part in coalescence events, thus strongly altering the
predictions of classical coalescent models. In particular, the
seed-bank coalescent ``does not come down from infinity,'' and the time
to the most recent common ancestor of a sample of size $n$ grows like
$\log\log n$.
This is in line with the empirical observation that seed-banks
drastically increase genetic variability in a population and indicates
how they may serve as a buffer against other evolutionary forces such
as genetic drift and selection.
\end{abstract}

%
\end{frontmatter}

\section{Introduction}
\label{sec:intro}

Seed-banks can play an important role in the population genetics of a
species, acting as a buffer against evolutionary forces such as random
genetic drift and selection as well as environmental variability (see,
e.g., \cite{ZT12,LJ11,Tellier} for an overview). Their presence typically
leads to significantly increased genetic variability respectively effective
population size (see, e.g., \cite{TL79,L90,N02,V04})
and could thus be considered as an important ``evolutionary force.'' In
particular,
classical mechanisms such as fixation and extinction of genes become
more complex: Genetic types can in principle disappear completely from
the \emph{active} population at a certain time while returning later
due to the germination of seeds or activation of dormant forms.

Seed-banks and dormant forms are known for many taxa. For example, they
have been suggested to play an important role in microbial evolution
\cite{LJ11,G14}, where certain bacterial endospores can remain
viable for (in principle arbitrarily) many generations.

Despite the many empirical studies
concerned with seed-bank effects, their mathematical modelling in
population genetics is still incomplete. Yet probabilistic models, and
in particular the Kingman coalescent (cf. \cite{K82}) and its relatives, have proven to
be very useful tools to understand basic principles of population
genetics and the interaction of evolutionary forces \cite{W09}.
Hence, it is natural to try to incorporate and investigate seed-banks
effects in such a probabilistic modelling framework, which is the aim
of this article. Before we present our model, we now briefly review
some of the mathematical work to date.

In 2001, Kaj, Krone and Lascoux \cite{KKL01} postulated and studied an
extension of the classical Wright--Fisher model from \cite{F30,W31} that includes
seed-banks effects. In their model, each generation consists of a fixed
amount of $N$ individuals. Each individual chooses its parent a random
amount of generations in the past and copies its genetic type. Here,
the number of generations that separates each parent and offspring is
understood as the time that the offspring spends as a seed or dormant
form. Formally, a parent is assigned to each individual in generation
$k$ by first sampling a random number $B$, which is assumed to be
independent and identically distributed for each individual, and then
choosing a parent uniformly among the $N$ individuals in generation
$k-B$ (note that the case $B \equiv1$ is just the classical
Wright--Fisher model).

Kaj, Krone and Lascoux then prove that if the seed-bank age
distribution $\mu$ of $B$ is restricted to finitely many generations
$\{
1, 2, \ldots, m\}$, where $m$ is independent of $N$, then the ancestral
process induced by the seed-bank model converges, after the usual
scaling of time by a factor $N$, to a time changed (delayed) Kingman
coalescent, where the coalescent rates are multiplied by $\beta^2:=
1/\E[B]^2$. An increase of the expected value of the seed-bank age
distribution thus further decelerates the coalescent, leading to an
increase in the effective population size. However, as observed by
\cite
{ZT12}, since the overall coalescent tree structure is retained, this
leaves the relative allele frequencies within a sample unchanged. In
this scenario
we thus speak of a ``weak'' seed-bank effect.

More generally, \cite{BG13} show that a sufficient condition for
convergence to the Kingman coalescent (with the same scaling and delay)
in this setup is that $E[B]<\infty$ with $B$ independent of $N$.

A different extension of the model by Kaj, Krone and Lascoux is
considered in~\cite{ZT12}, where the authors combine the seed-bank
model of \cite{KKL01} with fluctuations in population size. They point
out that substantial germ banks with small germination rates may buffer
or enhance the effect of the demography. This indicates that seed-bank
effects affect the interplay of evolutionary forces and can have
important consequences.

In this respect, \cite{BG13} show that \textit{strong} seed-bank effects
can lead to a behaviour which is very different from the Kingman
coalescent. In particular, if the seed-bank age distribution is
``heavy-tailed,'' say, $\mu(k)=L(k)k^{-\alpha}$, for $k\in\mathbb{N}$,
where $L$ is a slowly varying function, then if $\alpha<1$ the expected
time for the most recent common ancestor is infinite, and if $\alpha
<1/2$ two randomly sampled individuals do not have a common ancestor at
all with positive probability. Hence, this will not only delay, but
actually completely alter the effect of random genetic drift.

It can be argued that such extreme behaviour seems artificial
(although, as mentioned above, bacterial endospores may stay viable for
essentially unlimited numbers of generations). Instead, one can turn to
the case $B=B(N)$ and scale the seed-bank age distribution $\mu$ with
$N$ in order to understand its interplay with other evolutionary forces
on similar scales. For example, in \cite{G14} and \cite{B14} the
authors study a seed-bank model with $\mu=\mu(N)=(1-\varepsilon
)\delta
_1+\varepsilon\delta_{N^\beta}, \beta>0, \varepsilon\in(0,1)$. They
show that for $\beta<1/3$ the ancestral process converges, after
rescaling the time by the nonclassical factor $N^{1+2\beta}$, to the
Kingman coalescent, so that the expected time to the most recent common
ancestor is highly elevated in this scenario. However, in particular
since the above seed-bank Wright--Fisher model is highly non-Markovian,
the results in other parameter regimes, in particular $\beta=1$, are
still elusive.

To sum up, while there are mathematical results in the weak seed-bank
regime, it appears as if the ``right'' scaling regimes for stronger
seed-bank models, and the potentially new limiting coalescent
structures, have not yet been identified. This is in contrast to many
other population genetic models, where
the interplay of suitably scaled evolutionary forces (such as mutation,
genetic drift, selection and migration) often leads to elegant limiting
objects, such as the ancestral selection graph \cite{NK96}, or the
structured coalescent \cite{H97,N90}. A particular problem is the loss
of the Markov property in Wright--Fisher models with long genealogical ``jumps.''

In this paper, we thus propose a new Markovian Wright--Fisher type
seed-bank model that allows for a clear forward and backward scaling
limit interpretation. In particular, the forward limit in a bi-allelic
setup will consist of a pair of (S)DEs describing the allele frequency
process of our model, while the limiting genealogy, linked by a duality
result, is given by an apparently new coalescent structure which we
call \emph{seed-bank coalescent}. In fact, the seed-bank coalescent can
be thought of as a structured coalescent of a two island model in a
``weak migration regime,'' in which however coalescences are
\emph{completely blocked} in one island. Despite this simple description, the
seed-bank coalescent exhibits qualitatively altered genealogical
features, both in comparison to the Kingman coalescent and the
structured coalescent. In particular, we prove in Theorem~\ref
{thmm:notback} that the seed-bank coalescent ``does not come down from
infinity,'' and in Theorem~\ref{tmrca} that the expected time to the
most recent common ancestor of an $n$ sample is of
asymptotic order $\log\log n$
as $n$ gets large. Interestingly, this latter scale agrees with the one
for the Bolthausen--Sznitman coalescent identified by Goldschmidt and
Martin \cite{GM05}.\looseness=-1

Summarising, the seed-bank coalescent seems to be an interesting and
natural scaling limit for populations in the presence of a ``strong''
seed-bank effect. In contrast to previous genealogies incorporating
(weak) seed-bank effects, it is an entirely new coalescent structure
and not a time-change of Kingman's coalescent, capturing the essence of
seed-bank effects in many relevant situations. We conjecture that the
seed-bank coalescent is \emph{universal} in the sense that it is likely
to arise as the genealogy of other population models, such as Moran
models with a suitable seed-bank component.

The remainder of this paper is organised as follows.

In Section~\ref{sec:WF}, we discuss the Wright--Fisher model with a
seed-bank component that has a geometric age structure, and show that
its two bi-allelic frequency processes (for ``active'' individuals and
``seeds'') converge to a two-dimensional system of (S)DEs. We derive
their \emph{dual} process and employ this duality to compute the
fixation probabilities as $t \to\infty$ (in law) of the system.

In Section~\ref{sec:seedbank_coalescent}, we define the \emph{seed-bank
coalescent} corresponding to the previously derived dual block-counting
process and show how it describes the ancestry of the Wright--Fisher
geometric seed-bank model.

In Section~\ref{sec:properties}, we prove some interesting properties
of the seed-bank coalescent, such as ``\emph{not} coming down from
infinity'' and asymptotic bounds on the expected time to the most recent
common ancestor, which show that genealogical properties of a
population in the presence of strong seed-banks are altered qualitatively.

Sections~\ref{sec:seedbank_coalescent} and \ref{sec:properties} may be read essentially
independently of Section~\ref{sec:WF},
thus readers who are interested in the coalescent process only may
proceed directly to Section~\ref{sec:seedbank_coalescent}.

\section{The Wright--Fisher model with geometric seed-bank component}
\label{sec:WF}

\subsection{The forward model and its scaling limit}
\label{ssn:forward_model}

Consider a haploid population of fixed size $N$ reproducing in fixed
discrete generations $k=0,1,\ldots.$
Assume that individuals carry a genetic type from some type-space $E$
(we will later pay special attention to the bi-allelic setup, say $E=\{
a, A\}$, for the forward model).

Further, assume that the population also sustains a \emph{seed-bank} of
constant size $M=M(N)$, which consists of the dormant individuals. For
simplicity, we will frequently refer to the $N$ ``active'' individuals as
``plants'' and to the $M$ dormant individuals as ``seeds.''

Given $N, M\in\N$, let $\varepsilon\in[0,1]$ such that $\varepsilon N
\leq M$ and set $\delta:= \varepsilon N/M$,
and assume for convenience that $\varepsilon N=\delta M$ is a natural
number (otherwise replace it by $\lfloor\varepsilon N\rfloor$
everywhere). Let $[N]:=\{1, \ldots, N\}$ and $[N]_0:=[N] \cup\{0\}$.
The dynamics of our Wright--Fisher model with strong seed-bank
component are then as follows:
\begin{itemize}
\item The $N$ active individuals (plants) from generation 0 produce $
(1-\varepsilon) N$ active individuals in generation 1 by multinomial
sampling with equal weights.
\item Additionally, $\delta M= \varepsilon N$ uniformly (without
replacement) sampled seeds from the seed-bank of size $M$ in generation
0 ``germinate,'' that is, they turn into exactly one active individual in
generation 1 each, and leave the seed-bank.
\item The active individuals from generation 0 are thus replaced by
these $(1-\varepsilon) N + \delta M=N$ new active individuals, forming
the population of plants in the next generation 1.
\item Regarding the seed-bank, the $N$ active individuals from
generation 0 produce $\delta M=\varepsilon N$ seeds by multinomial
sampling with equal weights, filling the vacant slots of the seeds that
were activated.
\item The remaining $(1-\delta)M$ seeds from generation 0 remain
inactive and stay in the seed-bank (or, equivalently, produce exactly
one offspring each, replacing the parent).
\item Throughout reproduction, offspring and seeds copy/respectively
maintain the genetic type of the parent.
\end{itemize}
Thus, in generation 1, we have again $N$ active individuals and $M$
seeds. This probabilistic mechanism is then to be repeated
independently to produce generations $k=2,3,\ldots.$ Note that the
offspring distribution of active individuals (both for the number of
plants and for the number of seeds) is exchangeable within their
respective sub-population.
Further, one immediately sees that the time that a given seed stays in
the seed-bank before becoming active is geometric with success
parameter~$\delta$, while the probability a given plant produces a
dormant seed is $\varepsilon$.

\begin{figure}[b]

\includegraphics{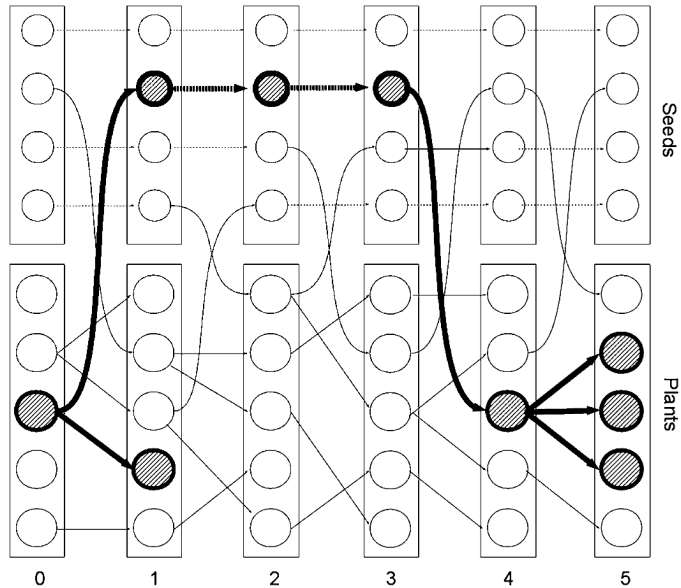}

\caption{A realisation of ancestral relationships in a
Wright--Fisher model with geometric seed-bank component.
Here, the genetic type of the third plant in generation 0 (highlighted
in grey) is lost after one generation,
but returns in generation four via the seed-bank, which acts as a
buffer against genetic drift and maintains genetic variability.}\label
{fig:WF_seedbank_model}
\end{figure}

\begin{defn}[(Wright--Fisher model with geometric seed-bank component)]
\label{def:forward}
Fix population size $N \in\N$, seed-bank size $M=M(N)$, genetic-type
space $E$
and $\delta, \varepsilon$ as before. Given initial
type configurations $\xi_0 \in E^{N}$ and $\eta_0 \in E^M$, denote by
\[
\xi_k:= \bigl(\xi_k(i) \bigr)_{i \in[N]},\qquad k \in\N,
\]
the random genetic-type configuration in $E^{N}$ of the plants in
generation $k$ (obtained from the above mechanism), and denote by
\[
\eta_k:= \bigl(\eta_k(j) \bigr)_{j \in[M]},\qquad k \in\N,
\]
correspondingly the genetic-type configuration of the seeds in $E^{M}$.
We call the discrete-time Markov chain $(\xi_k, \eta_k)_{k \in\N_0}$
with values in $E^{N} \times E^{M}$ the \emph{type configuration
process} of the \emph{Wright--Fisher model with geometric seed-bank component}.
See Figure \ref{fig:WF_seedbank_model} for a possible realisation of the model.
\end{defn}

We now specialise to the bi-allelic case $E=\{a,A\}$ and define the
frequency chains of $a$ alleles in the active population and in the
seed-bank. Define
%
\begin{equation}
\label{eq:frequency_chains} X_k^N:= \frac{1}N \sum
_{i\in[N]} {\mathbf1}_{\{\xi_k(i)=a\}}\quad \mbox{and}\quad
Y_k^M:= \frac{1}M \sum
_{j \in[M]} {\mathbf1}_{\{\eta_k(j)=a\}},\qquad k \in\N_0.\hspace*{-25pt}
\end{equation}
Both are discrete-time Markov chains taking values in
\begin{eqnarray*}
I^N&=& \biggl\{0, \frac{1}N, \frac{2}N, \ldots, 1
\biggr\} \subset[0,1] \quad\mbox{respectively}\\
 I^M&=& \biggl\{0,
\frac{1}M, \frac{2}M, \ldots, 1 \biggr\}\subset[0,1]. %
\end{eqnarray*}
Denote by $\P_{x,y}$ the distribution for which $(X^N, Y^M)$ starts in
$(x,y) \in I^N \times I^M$ $\P_{x,y}$-a.s., that is,
\[
\P_{x,y}( \cdot):= \P \bigl( \cdot | X^N_0 =
x, Y^M_0=y \bigr) \qquad\mbox {for } (x,y) \in
I^N \times I^M %
\]
(with analogous notation for the expectation, variance, etc.). The
corresponding time-homogeneous transition probabilities can now be
characterised.

\begin{prop}
\label{prop:transprob}
Let $c:=\varepsilon N = \delta M$ and assume $c \in[N]_0$. With the
above notation we have for $(x,y)$, respectively, $(\bar x, \bar y)
\in I^N \times I^M$,
\begin{eqnarray*}
p_{x,y}&:=&\P_{x,y} \bigl( X^N_1 = \bar
x, Y^M_1 = \bar y \bigr)
\\
&=& \sum_{i=0}^{c} \P_{x,y}( Z =
i)\P_{x, y}( U = \bar xN - i) \P _{x,y}(V = \bar yM + i),
\end{eqnarray*}
where $Z,U,V$ are independent under $\P_{x,y}$ with distributions
\[
\calL_{x,y}(Z) = \Hyp_{M,c,yM},\qquad \calL_{x,y}(U) =
\Bin_{N-c, x}, \qquad \calL _{x,y}(V) = \Bin_{c,x}.
\]
\begin{figure}

\includegraphics{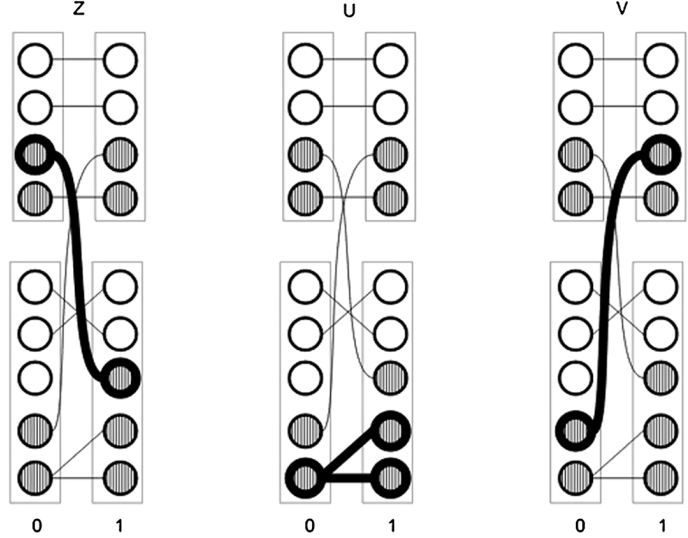}

\caption{In this figure $Z=1$, $U=2$ and $V=1$.}\label
{fig:RVrepresentation}\vspace*{-5pt}
\end{figure}
Here, $\Hyp_{M,c,yM}$ denotes the hypergeometric distribution with
parameters $M,c,  y\cdot M$ and $\Bin_{c,x}$ is the binomial distribution
with parameters $c$ and $x$.
\end{prop}

\begin{rem}
\label{rem:interpretation}
The random variables introduced in Proposition~\ref{prop:transprob}
have a simple interpretation, that is illustrated in Figure \ref{fig:RVrepresentation}:
\begin{longlist}[\textit{Z}]
\item[\textit{Z}] is the number of plants in generation 1, that are
offspring of a seed of type $a$ in generation 0. This corresponds to
the number of seeds of type $a$ that germinate/become active in the
next generation (noting that, in contrast to plants, the ``offspring'' of
a germinating seed is always precisely one plant and the seed vanishes).
\item[\textit{U}] is the number of plants in generation 1, that are
offspring of plants of type $a$ in generation 0.
\item[\textit{V}] is the number of seeds in generation 1, that are
produced by plants of type $a$ in generation 0.
\end{longlist}
\end{rem}

\begin{pf*}{Proof of Proposition~\ref{prop:transprob}}
With the interpretation of $Z,U$ and $V$ given in Remark~\ref
{rem:interpretation} their distributions are immediate as described in
the Definition~\ref{def:forward}. By construction, we then have $X^N_1
= \frac{U + Z}{N}$ and $Y^M_1 = {y}+\frac{V- Z}{M}$ and thus the
claim follows.
\end{pf*}

In many modelling scenarios in population genetics, parameters
describing evolutionary forces such as mutation, selection and
recombination are scaled in terms of the population size $N$ in order
to reveal a nontrivial limiting structure (see, e.g.,~\cite{E13} for
an overview). In our case, the interesting regime is reached
by letting $\varepsilon, \delta$ (and $M$) scale with $N$. More
precisely, assume that there exist $c,K\in(0,\infty)$ such that
%
\begin{equation}
\label{assumptions} \varepsilon=\varepsilon(N)=\frac{c}{N}\quad \mbox{and}\quad M=M(N)=
\frac
{N}{K}.
\end{equation}
Without loss of generality, $c\in[N]_0$ as $N\to\infty$. Under
assumption \eqref{assumptions}, the seed-bank age distribution is
geometric with parameter
%
\begin{equation}
\delta=\delta(N)=\frac{c}{M(N)}=\frac{cK}{N},
\end{equation}
and $c$ is the number of seeds that become active in each generation,
respectively, the number of individuals that move to the seed-bank. The
parameter $K$ determines the relative size of the seed-bank with
respect to the active population.

\begin{prop}
\label{prop:convergence}
Assume that \eqref{assumptions} holds. Consider test functions $f \in
C^{(3)}([0,1]^2)$. For any $(x,y) \in I^N \times I^M$, we define the
\emph{discrete generator} $A^N = A^N_{(\varepsilon, \delta, M)}$
of the frequency Markov chain $(X^N_k, Y^M_k)_{k \in\N}$ by 
%
\[
\label{eq:discrete_generator} A^N f(x,y):= N \E_{x,y} \bigl[ f
\bigl(X^N_1, Y^M_1 \bigr) - f(x,y)
\bigr].
\]
Then for all $(x,y)\in[0,1]^2$,
\[
\lim_{N\rightarrow\infty} A^Nf(x,y) = Af(x,y),
\]
where $A$ is defined by
\[
A f(x,y):= c (y-x)\frac{\partial f}{\partial x}(x,y) + c K(x-y)\frac
{\partial f}{\partial y} (x,y) +
\frac{1}{2} x(1-x)\frac{\partial
f^2}{\partial x^2}(x,y).
\]
\end{prop}

A proof can be found in the \hyperref[app]{Appendix}; see Proposition~\ref{prop:forgen}.
Since the state space of our frequency chain can be embedded in the
compact unit square $[0,1]^2$,
we get tightness and convergence on path-space easily by standard
argument (see, e.g., \cite{EK86}, Theorems~4.8.2 and~3.7.8) and
can identify the limit of our frequency chains as a pair of the
following (S)DEs.

\begin{corollary}[(Wright--Fisher diffusion with seed-bank component)]
\label{cor:sdes}
Under the conditions of Proposition~\ref{prop:convergence}, if
$X_0^N\to x$ \emph{a.s.} and $Y_0^M\to y$ \emph{a.s.}, we have that
\[
\bigl(X^N_{\lfloor Nt \rfloor}, Y^N_{\lfloor Nt \rfloor}
\bigr)_{t \geq0} \Rightarrow(X_t, Y_t)_{ t \geq0}
\]
on $D_{[0, \infty)}([0,1]^2)$ as $N\to\infty$, where $(X_t, Y_t)_{ t
\geq0}$ is a 2-dimensional diffusion solving
%
\begin{eqnarray}
\label{eq:system} d X_t & =& c(Y_t
-X_t)\,dt + \sqrt{X_t(1-X_t)}
\,dB_t,
\nonumber
\\[-8pt]
\\[-8pt]
\nonumber
d Y_t & =& cK(X_t -Y_t)\,dt,
\end{eqnarray}
where $(B_t)_{t\geq0}$ is standard Brownian motion.
\end{corollary}

The proof again follows from standard arguments, cf., for example,
\cite{EK86}, where in particular
Proposition~2.4 in Chapter~8 shows that the operator $A$ is indeed the
generator of a Markov process.

\begin{rem}
If we abandon the assumption $N=KM$, there are situations in which we
can still obtain meaningful scaling limits. If we assume
$N/M\rightarrow0$, and we rescale the generator as before by measuring
the time in units of size $N$, we obtain (cf. Proposition~\ref{prop:forgen})
\[
\lim_{N\rightarrow\infty} A^Nf(x,y) = c (y-x)
\frac{\partial
f}{\partial x}(x,y) +\frac{1}{2} x(1-x)\frac{\partial f^2}{\partial x^2}(x,y).
\]
This shows that the limiting process is purely one-dimensional, namely
the seed-bank frequency $Y_t$ is constantly equal to $y$, and the
process $(X_t)_{t\geq0}$ is a Wright--Fisher diffusion with migration
(with migration rate $c$ and reverting to the mean~$y$). The seed-bank,
which in this scaling regime is much larger than the active population,
thus acts as a reservoir with constant allele frequency $y$, with which
the plant population interacts.

The case $M/N\rightarrow0$ leads to a simpler limit: If we rescale the
generator by measuring the time in units of size $M$, we obtain
\[
\lim_{M\rightarrow\infty} A^Mf(x,y) = c (y-x)
\frac{\partial
f}{\partial y}(x,y)
\]
and constant frequency $X\equiv x$ in the plant population, which tells
us that if the seed-bank is of smaller order than the active
population, the genetic configuration of the seed-bank will converge to
the genetic configuration of the active population, in a deterministic way.
\end{rem}

The above results can be extended to more general genetic types spaces
$E$ in a standard way
using the theory of measure-valued respectively Fleming--Viot processes. This
will be treated elsewhere.
Before we investigate some properties of the limiting system, we first
derive its \emph{dual} process.

\subsection{The dual of the seed-bank frequency process}
\label{ssn:dual}

The classical Wright--Fisher diffusion is known to be dual to the block
counting process of the Kingman coalescent, and similar duality
relations hold for other models in population genetics. Such dual
processes are often extremely useful for the analysis of the underlying
system, and it is easy to see that our Wright--Fisher diffusion with
geometric seed-bank component also has a nice dual.

\begin{defn}
\label{Block}
We define the \emph{block-counting process of the seed-bank coalescent}
$(N_t,M_t)_{t \ge0}$ to be the continuous time Markov chain taking
values in $\N_0\times\N_0$ with
transitions
%
\begin{equation}
\label{eq:dual_rates} (n,m)\mapsto\cases{(n-1,m+1), &\quad$\mbox{at rate } cn,$ \vspace
*{2pt}
\cr
(n+1,m-1), & \quad$\mbox{at rate } cKm,$\vspace*{2pt}
\cr
(n-1,m), &\quad $\mbox{at rate } \pmatrix{n
\cr
2}$.}
\end{equation}
\end{defn}

Note that the three possible transitions correspond respectively to the
drift of the $X$-component, the drift of the $Y$-component, and the
diffusion part of the system~\eqref{eq:system}. This connection is
exploited in the following result.

Denote by $\P^{n,m}$ the distribution for which $(N_0,M_0)=(n,m)$ holds
$\P^{n,m}$-a.s., and denote the corresponding expected value by $\E
^{n,m}$. It is easy to see that, {eventually},
$N_t +M_t=1$ (as $t \to\infty$), $\P^{n,m}$-a.s. for all $n,m \in\N
_0$. We now show that $(N_t,M_t)_{t \ge0}$ is the \emph{moment dual}
of $(X_t,Y_t)_{t \ge0}$.

\begin{thmm}
\label{thmm:dual}
For every $(x,y)\in[0,1]^2 $, every $n,m\in\N_0$ and every $t\geq0$
%
\begin{equation}
\mathbb{E}_{x,y} \bigl[X_t^n
Y_t^m \bigr]=\mathbb{E}^{n,m}
\bigl[x^{N_t} y^{M_t} \bigr].
\end{equation}
\end{thmm}

\begin{pf}
Let $f(x,y;n,m):= x^ny^m $. Applying for fixed $n,m\in\N_0$ the
generator $A$ of $(X_t,Y_t)_{t \geq0}$ to $f$ acting as a function of
$x$ and $y$
gives
\begin{eqnarray*}
Af(x,y)&=&c(y-x)\frac{df}{dx}f(x,y)+\frac{1}{2}x(1-x)
\frac
{d^2f}{dx^2}f(x,y)\\
&&{}+cK(x-y)\frac{df}{dy}f(x,y)
\\
&=&c(y-x)nx^{n-1}y^m+\frac{1}{2}x(1-x)n(n-1)x^{n-2}y^m
\\
&&{}+cK(x-y)x^{n}my^{m-1}
\\
&=&cn\bigl(x^{n-1}y^{m+1}-x^ny^m\bigr)+
\pmatrix{n
\cr
2}\bigl(x^{n-1}y^{m}-x^ny^m
\bigr)
\\
&&{}+cKm\bigl(x^{n+1}y^{m-1}-x^ny^m
\bigr).
\end{eqnarray*}
Note that the right-hand side is precisely the generator of
$(N_t,M_t)_{t \ge0}$ applied to $f$ acting as a function of $n$ and
$m$, for fixed $x,y\in[0,1]$. Hence, the duality follows from standard
arguments; see, for example, \cite{JK14}, Proposition~1.2.
\end{pf}

\subsection{Long-term behaviour and fixation probabilities}
\label{ssn:fixation}

The long-term behaviour of our system \eqref{eq:system} is not obvious.
While a classical Wright--Fisher diffusion $(Z_t)_{t\geq0}$, given by
\[
dZ_t =\sqrt{Z_t(1-Z_t)}
\,dB_t, \qquad Z_0=z \in[0,1], %
\]
will get absorbed at the boundaries after finite time a.s. (in fact
with finite expectation), hitting 1 with probability $z$, this is more
involved for our frequency process in the presence of a strong
seed-bank. Nevertheless, one can still compute its fixation
probabilities as $t \to\infty$, at least in law. Obviously, $(0,0)$
and $(1,1)$ are absorbing states for the system \eqref{eq:system}. They
are also the only absorbing states, since absence of drift requires
$x=y$, and for the fluctuations to disappear, it is necessary to have
$x\in\{0,1\}$.

\begin{prop}
\label{prop:moment_limit}
All mixed moments of $(X_t,Y_t)_{t \ge0}$ solving \eqref{eq:system}
converge to the \emph{same} finite limit depending only on $x,y, K$.
More precisely, for each fixed $n,m\in\N$, we have
%
\begin{equation}
\label{eq:moment_value} \lim_{t \to\infty} \mathbb{E}_{x,y}
\bigl[X_t^{n}Y_t^{m}\bigr] =
\frac{y+xK}{1+K}.
\end{equation}
\end{prop}

\begin{pf}
Let $(N_t,M_t)_{t \ge0}$ be as in Definition~\ref{Block}, started in
$(n, m)\in\N_0\times\N_0$.
Let $T$ be the first time at which there is only one particle left in
the system $(N_t, M_t)_{t \ge0}$, that is,
\[
T:= \inf \{t >0\dvtx N_t+M_t=1 \}. %
\]
Note that for any finite initial configuration $(n,m)$, the stopping
time $T$ has finite expectation. Now, by Theorem~\ref{thmm:dual},
\begin{eqnarray*}
\lim_{t\rightarrow\infty}\E_{x,y} \bigl[X_t^{n}Y_t^{m}
\bigr] & =& \lim_{n\rightarrow\infty} \E^{n,m} \bigl[x^{N_t}y^{M_t}
\bigr]
\\
& =& \lim_{t\rightarrow\infty} \mathbb{E}^{n,m}
\bigl[x^{N_t}y^{M_t}\mid T\le t \bigr] \P^{n,m} (T\le t
)
\\
&&{} + \lim_{t\rightarrow\infty} \underbrace{\E^{n,m}
\bigl[x^{N_t}y^{M_t} |T>t \bigr]}_{\leq1}
\P^{n,m} (T>t )
\\
& =& \lim_{t\rightarrow\infty} \bigl(x \P^{n,m} (N_t=1,
T\leq t )+ y \P^{n,m} (M_t=1, T\leq t ) \bigr)
\\
& =& \lim_{t\rightarrow\infty} \bigl(x \P^{n,m} (N_t=1
)+ y \P ^{n,m} (M_t=1 ) \bigr)
\\
& = &\frac{xK}{1+K} + \frac{y}{1+K},
\end{eqnarray*}
where the last equality holds by convergence to the invariant
distribution of a single particle, jumping between the two states
``plant'' and ``seed'' at rate $c$ respectively $cK$, which is given by
$(K/(1+K), 1/(1+K))$ and independent of the choice of $n,m$.
\end{pf}

\begin{corollary}[(Fixation in law)]\label{cor:fix_law}
Given $c, K$, $(X_t, Y_t)$ converges in distribution as $t\to\infty$ to
a two-dimensional random variable $(X_\infty, Y_\infty)$, whose
distribution is given by
%
\begin{equation}
\label{eq:momentconvergence} \mathcal{L}_{(x,y)} ( X_\infty,
Y_\infty ) = \frac{y+xK}{1+K} \delta_{(1,1)} +
\frac{1+(1-x)K-y}{1+K} \delta_{(0,0)}.
\end{equation}
\end{corollary}

Note that this is in line with the classical results for the
Wright--Fisher diffusion: As $K \to\infty$ (i.e., the seed-bank
becomes small compared to the plant population), the fixation
probability of
$a$ alleles approaches $x$. Further, if $K$ becomes small (so that the
seed-bank population dominates the plant population), the fixation
probability is governed by the initial fraction $y$ of $a$-alleles in
the seed-bank.

\begin{pf*}{Proof of Corollary \ref{cor:fix_law}}
It is easy to see that the only two-dimen\-sional distribution on
$[0,1]^2$, for which all moments are constant equal to $\frac
{xK+y}{1+K}$, is given by
\[
\frac{y+xK}{1+K} \delta_{(1,1)} + \frac{1+(1-x)K-y}{1+K}
\delta_{(0,0)}. %
\]
Indeed, uniqueness follows from the moment problem, which is uniquely
solvable on $[0,1]^2$. Convergence in law follows from convergence of
all moments due to Theorem~3.3.1 in \cite{EK86} and the
Stone--Weierstra\ss\ theorem.
\end{pf*}

\begin{rem}[(Almost sure fixation)]
Observing that $(KX_t+Y_t)_{t\geq0}$ is a bounded martingale, and
given the shape of the limiting law \eqref{eq:momentconvergence}, one
can also get almost sure convergence of $(X_t, Y_t)$ to $(X_\infty,
Y_\infty)$ as $t\to\infty$. However, as we will see later, fixation
will not happen in finite time, since the
block-counting process $(N_t,M_t)_{t \ge0}$, started from an infinite
initial state, \emph{does not come down from infinity} (see
Section~\ref
{sec:properties}), which means that the whole (infinite) population
does not have a most-recent common ancestor. Thus, in finite time,
initial genetic variability should never be completely lost. We expect
that with some extra work, this intuitive reasoning could be made
rigorous in an almost sure sense with the help of a ``look-down
construction,'' and will be treated in future work. The fact that
fixation does not occur in finite time can also be understood from
\eqref
{eq:system}, where we can compare the seed-component $(Y_t)_{t\geq0}$
to the solution of the deterministic equation
\[
dy_t=-cKy_t\,dt,
\]
corresponding to a situation where the drift towards 0 is maximal [or
to $dy_t=cK(1-y_t)\,dt$ where the drift towards 1 is maximal]. Since
$(y_t)_{t\geq0}$ does not reach~0 in finite time if $y_0>0$, neither
does $(Y_t)_{t\geq0}$.
\end{rem}

\section{The seed-bank coalescent}
\label{sec:seedbank_coalescent}

\subsection{Definition and genealogical interpretation}
In view of the form of the block counting process, it is now easy to
guess the stochastic process describing the limiting gene genealogy of
a sample taken from the Wright--Fisher model with seed-bank component.
Indeed, for $k \ge1$, let $\mathcal{P}_k$ be the set of partitions of
$[k]$. For $\pi\in\mathcal{P}_k$ let $|\pi|$ be the number of blocks
of the partition $\pi$. We define the space of \emph{marked} partitions
to be
\[
\mathcal{P}^{\{p,s\}}_k= \bigl\{ (\zeta, \vec{u}) \mid\zeta\in
\mathcal {P}_k, \vec{u} \in\{s,p\}^{|\zeta|} \bigr\}. %
\]
This enables us to attach to each partition block a flag which can be
either ``plant'' or ``seed'' ($p$ or $s$), so that we can trace whether an
ancestral line is currently in the active or dormant part of the
population. For example, for $k=5$, an element $\pi$ of $\mathcal
{P}^{\{
p,s\}}_k$ is the marked partition $\pi= \{\{1,3\}^p\{2\}^s\{4,5\}
^p \}$.

\begin{figure}[b]

\includegraphics{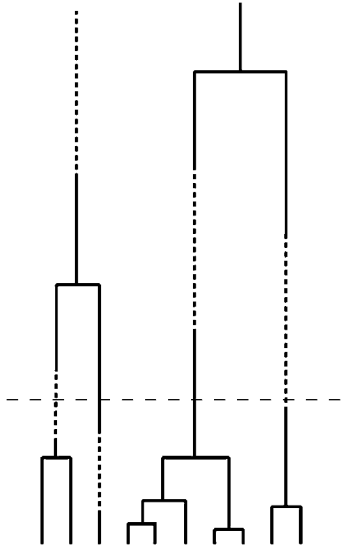}

\caption{A possible realisation of the standard 10-seed-bank
coalescent. Dotted lines indicate
``inactive lineages'' (carrying an $s$-flag, which are prohibited from merging).
At the time marked with the dotted horizontal line the process is in
the state $\{\{1,2\}^s\{3\}^p\{4,5,6,7,8\}^p\{9,10\}^s\}$.}\label{fig3}
\end{figure}

Consider two marked partitions $\pi, \pi^\prime\in\mathcal{P}_k^{\{
p,s\}}$, we say $\pi\succ\pi^\prime$ if $\pi^\prime$ can be
constructed by merging exactly 2 blocks of $\pi$ carrying the $p$-flag,
and the resulting block in ${\pi}^\prime$ obtained from the merging
both again carries a $p$-flag. For example,
\[
\bigl\{\{1,3\}^p\{2\}^s\{4,5\}^p \bigr\}
\succ \bigl\{\{1,3,4,5\}^p\{2\} ^s \bigr\}. %
\]
We use the notation ${\pi}\Join{\pi}^\prime$ if ${\pi
}^\prime$ can be constructed by changing the flag of precisely one
block of $\pi$, for example,
\[
\bigl\{\{1,3\}^p\{2\}^s\{4,5\}^p\bigr\}\Join
\bigl\{\{1,3\}^s\{2\}^s\{4,5\}^p \bigr\}.
\]

\begin{defn}[(The seed-bank $k$-coalescent)]
\label{defn:k_seed-bank_coalescent}
For $k \ge2$ and $c,K \in(0,\infty)$, we define the \emph{seed-bank
$k$-coalescent} $(\Pi^{(k)}_t)_{t \ge0}$ with seed-bank intensity $c$
and relative seed-bank size $1/K$ to be the continuous time Markov
chain with values in $\mathcal{P}_k^{\{p,s\}}$, characterised by the
following transitions:
%
\begin{equation}\qquad
\label{eq:coalescent_transitions}
{\pi} \mapsto{\pi}^\prime\mbox{ at rate }
\cases{1, &\quad $\mbox{if } {\pi}\succ{\pi}^\prime$,\vspace*{2pt}
\cr
c, &\quad $\mbox{if } {\pi}\Join{\pi}^\prime\mbox{ and one $p$
is replaced by one $s$}$,\vspace*{2pt}
\cr
cK, &\quad $\mbox{if } {\pi}\Join{\pi}^\prime\mbox{ and one $s$ is replaced by one $p$}$. }
\end{equation}
If $c=K=1$, we speak of the \emph{standard seed-bank $k$-coalescent}.
See Figure \ref{fig3} for a possible realisation.
\end{defn}

Comparing \eqref{eq:coalescent_transitions} to \eqref{eq:dual_rates},
it becomes evident that $(N_t, M_t)$ introduced in Definition~\ref
{Block} is indeed the block-counting process of the seed-bank coalescent.

\begin{defn}[(The seed-bank coalescent)]
\label{defn:projective_limit}
We may define the \emph{seed-bank coalescent}, $(\Pi_t)_{t \ge
0}=(\Pi
^{(\infty)}_t)_{t \ge0}$ with seed-bank intensity $c$ and relative
seed-bank size $1/K$ as the
unique Markov process distributed according to the projective limit as
$k$ goes to infinity of the laws of the seed-bank $k$-coalescents (with
seed-bank intensity $c$ and relative seed-bank size $1/K$).
In analogy to Definition~\ref{defn:k_seed-bank_coalescent}, we call the
case of $c=K=1$ the \emph{standard seed-bank coalescent}.
\end{defn}

\begin{rem}
\label{rem:well-defined}
Note that the seed-bank coalescent is a well-defined object. Indeed,
for the projective limiting procedure to make sense, we need to show
\emph{consistency} and then apply the Kolmogorov extension theorem.
This can be roughly sketched as follows. Define the process
$(\overrightarrow{\Pi}{}^{(k)}_t)_{t \ge0}$ as the projection of $(\Pi
^{(k+1)}_t)_{t \ge0}$, the $k+1$ seed-bank coalescent, to the space
$\mathcal{P}_k^{\{p,s\}}$. Mergers and flag-flips involving the
singleton $\{k+1\}$ are only visible in $(\Pi^{(k+1)}_t)_{t\ge0}$, but
do not affect $(\overrightarrow{\Pi}^{(k)}_t)_{t\ge0}$. Indeed, by the
Markov-property, a change involving the singleton $\{k+1\}$ does not
affect any of the other transitions.
Hence, if $\overrightarrow{\Pi}{}^{(k)}_0=\Pi^{(k)}_0$, then
\[
\bigl(\overrightarrow{\Pi}{}^{(k)}_t\bigr)_{t\ge0}=
\bigl(\Pi^{(k)}_t\bigr)_{t\ge0}
\]
holds in distribution.
By the Kolmogorov extension theorem, the projective limit exists and is unique.
\end{rem}

Note that it is obvious that the distribution of the block counting
process of the seed-bank coalescent, counting the number of blocks
carrying the $p$ and $s$-flags, respectively, agrees with the
distribution the process $(N_t, M_t)_{t \ge0}$ from Definition~\ref
{Block} (with suitable initial conditions).

Further, it is not hard to see that the seed-bank coalescent appears as
the limiting genealogy of a sample taken from the Wright--Fisher model
with geometric seed-bank component in the same way as the Kingman
coalescent describes the limiting genealogy of a sample taken from the
classical Wright--Fisher model (here, we merely sketch a proof, which
is entirely standard).

Indeed, consider the genealogy of a sample of $k\ll N$ individuals,
sampled from present generation $0$. We proceed backward in time,
keeping track in each generation of the ancestors of the original
sample among the active individuals (plants) and among the seeds. To
this end, denote by $\Pi_i^{(N,k)}\in\mathcal P^{\{p,s\}}_k$ the
configuration of the genealogy at generation $-i$, where two
individuals belong to the same block of the partition $\Pi_i^{(N,k)}$
if and only if their ancestral lines have met until generation $-i$,
which means that all individuals of a block have exactly one common
ancestor in this generation, and the flag $s$ or $p$ indicates whether
said ancestor is a plant or a seed in generation $-i$. According to our
forward in time population model, there are the following possible
transitions from one generation to the previous one of this process:
\begin{itemize}
\item One (or several) plants become seeds in the previous generation.
\item One (or several) seeds become plants in the previous generations.
\item Two (or more) individuals have the same ancestor in the previous
generation (which by construction is necessarily a plant), meaning that
their ancestral lines merge.
\item Any possible combination of these three events.
\end{itemize}

It turns out that only three of the possible transitions play a role in
the limit as $N\to\infty$, whereas the others have a probability that
is of smaller order.

\begin{prop}\label{prop:trans_prob_gene}
In the setting of Proposition~\ref{prop:convergence}, additionally
assume that $\Pi^{(N,k)}_0=\{\{1\}^p,\ldots,\{k\}^p\}$, $\P$-a.s. for
some fixed $k\in\N$. Then for $\pi, \pi'\in\mathcal P^{\{p,s\}}_k$,
%
\begin{eqnarray}
&&\P \bigl(\Pi_{i+1}^{(N,k)}=\pi' |
\Pi_{i}^{(N,k)}=\pi \bigr)
\nonumber
\\[-8pt]
\\[-8pt]
\nonumber
&&\qquad
=\cases{\displaystyle\frac{1}{N}+O
\bigl(N^{-2}\bigr), &\quad $\mbox{if }\pi\succ\pi'$,
\vspace*{2pt}
\cr
\displaystyle\frac{c}{N}+O\bigl(N^{-2}\bigr), &\quad $\mbox{if } {\pi}\Join{\pi}^\prime \mbox{ and a $p$ is replaced by an
$s$,}$\vspace*{2pt}
\cr
\displaystyle\frac{cK}{N}+O\bigl(N^{-2}\bigr), &\quad $ \mbox{if } {\pi}\Join{\pi }^\prime \mbox{ and an $s$ is replaced by
a $p$,}$\vspace*{2pt}
\cr
O\bigl(N^{-2}\bigr),&\quad $\mbox{otherwise}$}
\end{eqnarray}
for all $i \in\N_0$.
\end{prop}

\begin{pf}
According to the definition of the forward in time population model,
exactly $c$ out of the $N$ plants become seeds, and exactly $c$ out of
the $M=N/K$ seeds become plants. Thus, whenever the current state $\Pi
^{(N,k)}_i$ of the genealogical process contains at least one
$p$-block, then the probability that a \emph{given} $p$-block changes
flag to $s$ at the next time step is equal to $\frac{c}{N}$. If there
is at least one $s$-block, then the probability that any given
$s$-block changes flag to $p$ is given by $\frac{cK}{N}$, and the
probability that a given $p$-block chooses a \emph{fixed} plant
ancestor is equal to $ (1-\frac{c}{N} )\frac{1}{N}$ (where
$1-c/N$ is the probability that the ancestor of the block in question
is a plant, and $1/N$ is the probability to choose one particular plant
among the $N$).

From this, we conclude that the probability of a coalescence of two
given $p$-blocks in the next step is
\[
\P (\mbox{two given $p$-blocks merge} )= \biggl(1-\frac
{c}{N}
\biggr)^2\frac{1}{N}.
\]
Since we start with $k$ blocks, and the blocks move independently, the
probability that two or more blocks change flag at the same time is of
order at most $N^{-2}$. Similarly, the probability of any combination
of merger or block-flip events other than single blocks flipping or
binary mergers is of order $N^{-2}$ or smaller, since the number of
possible events (coalescence or change of flag) involving at most $k$
blocks is bounded by a constant depending on $k$ but not on $N$.
\end{pf}

\begin{corollary}
For any $k\in\N$, under the assumptions of Proposition~\ref
{prop:convergence}, $(\Pi^{(N,k)}_{\lfloor Nt\rfloor})_{t\geq0}$
converges weakly as $N\to\infty$ to the seed-bank coalescent $(\Pi
^{(k)}_t)_{t\geq0}$ started with $k$ plants.
\end{corollary}

\begin{pf}
From Proposition~\ref{prop:trans_prob_gene}, it is easy to see that the
generator of $(\Pi^{(N,k)}_{\lfloor Nt\rfloor})$ converges to the
generator of $(\Pi^{(k)}_t)$, which is defined via the rates given in
\eqref{eq:coalescent_transitions}. Then standard results
(see Theorem~3.7.8 in \cite{EK86}) yield weak convergence of the process.
\end{pf}

\subsection{Related coalescent models}
\label{ssect:relmod}

\textit{The structured coalescent}.
The seed-bank coalescent is reminiscent of the \emph{structured
coalescent} arising from a two-island population model (see, e.g.,
\cite{W31,T88,N90,H94,H97}). Indeed, consider two Wright--Fisher type
(sub-)populations of fixed relative size evolving on separate
``islands'', where individuals (resp., ancestral lineages) may migrate
between the two locations with a rate of order of the reciprocal of the
total population size (the so-called ``weak migration regime''). Since
offspring are placed on the same island as their parent, mergers
between two ancestral lineages are only allowed if both are currently
in the same island.
This setup again gives rise to a coalescent process defined on ``marked
partitions,'' with the marks indicating the location of the ancestral
lines among the two islands. Coalescences are only allowed for lines
carrying the same mark at the same time, and marks are switched
according to the scaled migration rates. See \cite{W09} for an overview.

In our Wright--Fisher model with geometric seed-bank component, we
consider a similar ``migration'' regime
between the two sub-populations, in our case called ``plants'' and
``seeds.'' However, in the resulting seed-bank coalescent, coalescences
can only happen while in the plant-population. This asymmetry leads to
a behaviour that is qualitatively different to the usual two-island
scenario (e.g., with respect to the time to the most recent common
ancestor, whose expectation is always finite for the structured
coalescent, even if the sample size goes to infinity).

\textit{The coalescent with freeze}.
Another related model is the \emph{coalescent with freeze} (see \cite
{DGP07}), where blocks become completely inactive at some rate. This
model is different from ours because once a block has become inactive,
it cannot be activated again. Hence, it cannot coalesce at all, which
clearly leads to a different long-time behaviour. In particular, one
will not expect to see a most recent common ancestor in such a coalescent.

\section{Properties of the seed-bank coalescent}
\label{sec:properties}

\subsection{Coming down from infinity}
\label{ssn:coming_down}

The notion of \emph{coming down from infinity} was discussed by Pitman
\cite{P99} and Schweinsberg \cite{S00}. They say that an exchangeable
coalescent process \emph{comes down from infinity} if the corresponding
block counting process (of an infinite sample) has finitely many blocks
immediately after time~0 (i.e., the number of blocks is finite almost
surely for each $t>0$). Further, the coalescent is said to \emph{stay
infinite} if the number of blocks is infinite a.s. for all $t \ge
0$. Schweinsberg also gives a necessary and sufficient criterion for
so-called ``Lambda-coalescents'' to come down from infinity. In
particular, the Kingman coalescent does come down from infinity.
However, note that the seed-bank coalescent does not belong to the
class of Lambda-coalescents, so that Schweinsberg's result does not
immediately apply.
For an overview of the properties of general exchangeable coalescent
processes, see, for example, \cite{B09}.

\begin{thmm}
\label{thmm:notback}
The seed-bank coalescent does \emph{not} come down from infinity. In
fact, its block-counting process $(N_t,M_t)_{t \ge0}$
stays infinite for every $t \geq0$, $\P$-a.s. To be precise, for each
starting configuration $(n, m)$ where $n+m$ is (countably) infinite,
\[
\P \bigl(\forall t\geq0\dvtx M^{(n,m)}_t= \infty \bigr) = 1.
\]
\end{thmm}

The proof of this theorem is based on a coupling with a dominated
simplified \emph{coloured seed-bank coalescent} process introduced
below. In essence, the coloured seed-bank coalescent behaves like the
normal seed-bank coalescent, except we mark the individuals with a
colour to indicate whether they have (entered and) left the seed-bank
at least once. This will be useful in order to obtain a process where
the number of plant-blocks is nonincreasing. We will then prove that
even if we consider only those individuals that have never made a
transition from seed to plant (but possibly from plant to seed), the
corresponding block-counting process will stay infinite. This will be
achieved by proving that infinitely many particles enter the seed-bank
before any positive time. Since they subsequently leave the seed-bank
at a linear rate, this will take an infinite amount of time.

\begin{defn}[(A coloured seed-bank coalescent)]
In analogy to the construction of the seed-bank coalescent, we first
define the set of \emph{coloured}, marked partitions as
\begin{eqnarray*}
\mathcal P^{\{p,s\}\times\{w,b\}}_k &:=& \bigl\{(\pi, \vec{u}, \vec{v}) \mid(
\pi, \vec{u}) \in\mathcal P^{\{p,s\}}_k, \vec{v} \in\{w,b\}
^k \bigr\},\qquad k \in\N,
\\
\mathcal P^{\{p,s\}\times\{w,b\}} &:= &\bigl\{(\pi, \vec{u}, \vec{v}) \mid(\pi, \vec{u})
\in\mathcal P^{\{p,s\}}, \vec{v} \in\{w,b\} ^{\N
} \bigr\}.
\end{eqnarray*}
It corresponds to the marked partitions introduced earlier, where now
each element of $[k]$, respectively, $\N$, has an additional flag
indicating its colour: $w$ for \emph{white} and $b$ for \emph{blue}. We
write $\pi\succ\pi'$, if $\pi'$ can be constructed from $\pi$ by
merging two blocks with a $p$-flag in $\pi$ that result into a block
with a $p$-flag in $\pi'$, while each individual retains its colour. It
is important to note that the $p$- or $s$-flags are assigned to \emph
{blocks}, the colour-flags to \emph{individuals}, that is, elements of
$[k]$, respectively, $\N$.
We use $\pi\ltimes\pi'$, to denote that $\pi'$ results from $\pi$ by
changing the flag of a block from $p$ to $s$ and leaving the colours of
all individuals unchanged and $\pi\rtimes\pi'$, if $\pi'$ is obtained
from $\pi$, by changing the flag of a block from $s$ to $p$ \emph{and
colouring all the individuals in this block blue}, that is, setting
their individual flags to $b$. In other words, after leaving the
seed-bank, individuals are always coloured blue.

For $k \in\N$ and $c,K \in(0,\infty)$ we now define the \emph
{coloured seed-bank $k$-coalescent with seed-bank intensity c and
seed-bank size $1/K$}, denoted by $(\underline\Pi_t)_{t \geq0}$, as the
continuous time Markov chain with values in $\mathcal P^{\{p,s\}\times
\{
w,b\}}_k$ and transition rates given by
%
\begin{equation}
\pi\mapsto\pi' \mbox{ at rate }\cases{1,&\quad $\mbox{if } \pi\succ
\pi '$,\vspace*{2pt}
\cr
c, &\quad $\mbox{if } \pi\ltimes
\pi'$, \vspace*{2pt}
\cr
cK, &\quad $\mbox{if } \pi\rtimes
\pi'$.}
\end{equation}
The \emph{coloured seed-bank coalescent with seed-bank intensity c and
seed-bank size $1/K$} is then the unique Markov process on $\mathcal P^{\{
p,s\}\times\{w,b\}}$ given by the projective limit of the distributions
of the $k$-coloured seed-bank coalescents, as $k$ goes to infinity.
\end{defn}

\begin{rem}
\label{rem:csbc}
1. Note that the coloured seed-bank coalescent is well-defined.
Since the colour of an individual only depends on its own path and does
not depend on the colour of other individuals (not even those that
belong to the same block), the consistency of the laws of the
$k$-coloured seed-bank coalescents boils down to~the consistency of the
seed-bank $k$-coalescents discussed in Remark~\ref{rem:well-defined}.
In much the same way, we then obtain the existence and uniqueness of
the coloured seed-bank coalescent from Kolmogorov's Extension theorem.\vspace*{-6pt}
\begin{longlist}[2.]
\item[2.] The normal seed-bank ($k$-)coalescent can be obtained from the
coloured seed-bank ($k$-)coalescent by omitting the flags indicating
the colouring of the individuals. However, if we only consider those
blocks containing \emph{at least} one white individual, we obtain a
coalescent similar to the seed-bank coalescent, where lineages are
discarded once they leave the seed-bank.
\end{longlist}
\end{rem}

For $t \geq0$ define $\underline N_t$ to be the number of \emph{white
plants} and $\underline M_t$ the number of \emph{white seeds} in
$\underline\Pi_t$. We will use a superscript $(n,m)$ to denote the
processes started with $n$ plants and $m$ seeds $\P$-a.s., where $n,m =
\infty$ means we start with a countably infinite number of plants,
respectively, seeds. We will always start in a configuration were all
individual labels are set to $w$, that is, with only white particles.
Note that our construction is such that $(\underline{N}_t)_{t \ge0}$
is nonincreasing.

\begin{prop}
\label{prop:csbc}
For any $n, m \in\N\cup\{\infty\}$, the processes
$(N^{(n,m)}_t,\break M^{(n,m)}_t)_{t \geq0}$ and $(\underline
{N}^{(n,m)}_t,\underline{M}^{(n,m)}_t)_{t \geq0}$ can be coupled such that
\[
\P \bigl( \forall t \geq0\dvtx N^{(n,m)}_t \geq\underline
{N}^{(n,m)}_t \mbox{and } M^{(n,m)}_t
\geq\underline {M}^{(n,m)}_t \bigr) =1.
\]
\end{prop}

\begin{pf}
This result is immediate if we consider the coupling through the
coloured seed-bank coalescent and the remarks in Remark~\ref{rem:csbc}.
\end{pf}

\begin{pf*}{Proof of Theorem~\ref{thmm:notback}}
Proposition~\ref{prop:csbc} implies that it suffices to prove the
statement for $(\underline{M}_t)_{t \geq0}$ instead of $(M_t)_{t \geq
0}$. In addition, we will only have to consider the case of $m=0$,
since starting with more (possibly infinitely many) seeds will only
contribute towards our desired result.

For $n \in\N\cup\{\infty\}$, let
\[
\tau_j^n:=\inf\bigl\{t\geq0\dvtx\underline{N}^{(n,0)}_t=j
\bigr\},\qquad 1\leq j\leq n-1, j<\infty
\]
be the first time that the number of active blocks of an $n$-sample
reaches $k$. Note that $(\underline{N}_t)_{t \geq0}$ behaves like the
block-counting process of a Kingman coalescent where in addition to the
coalescence events, particles may ``disappear'' at a rate proportional
to the number of particles alive. Since the corresponding values for a
Kingman coalescent are finite $\P$-a.s., it is easy to see that the
$\tau^n_j$ are, too.
Clearly, for any $n$, $\tau^n_j-\tau^n_{j-1}$ has an exponential
distribution with parameter
\[
\lambda_j:=\pmatrix{j
\cr
2}+cj.
\]
At each time of a transition $\tau_j^n$, we distinguish between two
events: \emph{coalescence} and \emph{deactivation} of an active block,
where by deactivation we mean a transition of $(\underline{N}^n_t,
\underline{M}^n_t)_{t \geq0}$ of type $(j+1,l)\mapsto(j, l+1)$ (for
suitable $l \in[n]$), that is, the transition of a plant to a seed.

Then
%
\begin{equation}
\label{eq:deactivation-probability} \P \bigl( \mbox{deactivation at }\tau^n_{j-1}
\bigr) =\frac{cj}{{j\choose2}+cj} =\frac{2c}{j+2c-1},
\end{equation}
independently of the number of inactive blocks.
Thus,
\[
X^n_j:={\mathbf1}_{ \{\mathrm{deactivation\ at\ }\tau^n_{j-1} \}},\qquad j=2,\ldots,n, j<
\infty,
\]
are independent Bernoulli random variables with respective parameters
$2c/\break (j+2c-1), j=2,\ldots,n$. Note that $X_j^n$ depends on $j$, but the
random variable is independent of the random variable $\tau_{j-1}$ due
to the memorylessness of the exponential distribution.
Now define $A_t^n $ as the (random) number of deactivations up to time
$t\geq0$ that is, for $n\in\N\cup\{\infty\}$,
%
\begin{equation}
\label{def:An} A_t^n:=\sum
_{j=2}^n X^n_j {
\mathbf1}_{\{\tau^n_{j-1}<t\}}.
\end{equation}
For $n\in\N$, since $\lambda_j\geq{j\choose 2}$, it follows from
a comparison with the block counting process of the Kingman coalescent,
denoted by $(|\tilde\Pi_t^n|)_{t \ge0}$ if started in $n$ blocks, that
for all $t\geq0$,
\begin{eqnarray*}
\lim_{n\to\infty}\P \bigl(\tau^n_{\lfloor\log n-1\rfloor}\leq t
\bigr) & \geq&\lim_{n\to\infty} \P \bigl( \bigl|\tilde\Pi_t^n\bigr|
\le\lfloor\log n -1\rfloor \bigr)
\\
& \geq&\lim_{n\to\infty}\P \bigl( |\tilde\Pi_t| \le\log n -1
\bigr) =1,
\end{eqnarray*}
where the last equality follows from the fact that the
Kingman coalescent $(\tilde\Pi_t)_{t \geq0}$ comes down from infinity,
cf. \cite{S00,P99}.
For $t \ge0$,
%
\begin{eqnarray}
\label{eq:timelogn} \P \Biggl(A_t^n\geq\sum
_{j=\log n}^n X^n_j \Biggr) &
\ge&\P \Biggl({\mathbf 1}_{\{
\tau^n_{\log n-1}<t\}}\sum_{j=\log n}^n
X^n_j\geq\sum_{j=\log n}^n
X^n_j \Biggr)
\\
& \ge&\P \bigl(\tau^n_{\log n-1}<t \bigr)
\end{eqnarray}
and hence, by \eqref{eq:timelogn}
%
\begin{equation}
\label{eq:compA} \lim_{n \rightarrow\infty}\P \Biggl(A_t^n
\geq\sum_{j=\log
n}^nX^n_j
\Biggr) = 1.
\end{equation}
Note that due to \eqref{eq:deactivation-probability},
%
\begin{equation}
\label{eq:exp_An} \E \Biggl[\sum_{j=\log n}^n
X^n_j \Biggr]=\sum_{j=\log n}^n
\frac
{2c}{j+2c-1}=2c(\log n-\log\log n)+R(c,n),
\end{equation}
where $R(c,n)$ converges to a finite value depending on the seed-bank
intensity $c$ as $n\to\infty$. Since the $X^n_j$ are independent
Bernoulli random variables, we obtain for the variance
%
\begin{eqnarray}
\label{eq:var_An} \mathbb{V} \Biggl[\sum_{j=\log n}^nX^n_j
\Biggr]&=&\sum_{j=\log
n}^n\mathbb {V}
\bigl[X^n_j \bigr]=\sum_{j=\log n}^n
\frac{2c}{j+2c-1} \biggl(1-\frac
{2c}{j+2c-1} \biggr)
\nonumber
\\[-8pt]
\\[-8pt]
\nonumber
&\leq& 2c\log n\qquad\mbox{as } n\to\infty.
\end{eqnarray}
For any $\varepsilon>0$, we can choose $n$ large enough such that, $\E
[\sum_{j=\log n}^n X_k] \geq(2 c-\varepsilon) \log n$ holds, which yields
%
\begin{eqnarray}
\label{eq:chebyshev} \P \Biggl(\sum_{j=\log n}^n
X^n_j<c\log n \Biggr) & \leq&\P \Biggl(\sum
_{j=\log
n}^n X^n_j-\E \Biggl[
\sum_{j=\log n}^n X^n_j
\Biggr] < - (c-\varepsilon )\log n \Biggr)
\nonumber
\\
& \leq&\P \Biggl( \Biggl|\sum_{j=\log n}^n
X^n_j-\E \Biggl[\sum_{j=\log n}^n
X_j^n \Biggr] \Biggr|>(c-\varepsilon)\log n \Biggr)
\\
& \leq&\frac{2c}{(c-\varepsilon)^2\log n}\nonumber,
\end{eqnarray}
by Chebyshev's inequality.
In particular, for any $\kappa\in\mathbb{N}$,
\[
\lim_{n\rightarrow\infty}\P \Biggl(\sum_{j=\log n}^n
X^n_j<\kappa \Biggr)= 0,
\]
and together with \eqref{eq:compA} we obtain for any $t>0$
%
\begin{equation}
\label{eq:inf-seeds} \lim_{n\to\infty}\P\bigl(A^n_t
< \kappa\bigr)=0.
\end{equation}
Since the $(A^n_t)_{t\geq0}$ are coupled by construction for any $n
\in\N\cup\{\infty\}$, we know in particular that $\P(A^{\infty}_t
<\kappa) \leq\P(A^{n}_t <\kappa)$, for any $n \in\N, t\geq0,
\kappa
\geq0$ and therefore $\P(A^{\infty}_t < \kappa)=0$, which yields
%
\begin{equation}
\label{resultA} \forall t\geq0\qquad \P\bigl(A^{\infty}_t =
\infty\bigr) = 1.
\end{equation}
Since in addition, $(A^{\infty}_t)_{t\geq0}$ is nondecreasing in $t$,
we can even conclude
%
\begin{equation}
\label{resultB} \P\bigl(\forall t\geq0\dvtx A^{\infty}_t =
\infty\bigr) = 1.
\end{equation}
Thus, we have proven that, for any time $t \geq0$, there have been an
infinite amount of movements to the seed-bank $\P$-a.s. Now we are left
to show that this also implies the presence of an infinite amount of
lineages in the seed-bank, that is, that a sufficiently large
proportion is saved from moving back to the plants where it would be
``instantaneously'' reduced to a finite number by the coalescence mechanism.

Define $\mathcal B_t$ to be the blocks of a partition that visited the
seed-bank at some point before a fixed time $t\geq0$ and were visible
in the ``white'' seed-bank coalescent, that is,
\begin{eqnarray*}
&&\mathcal B_t:=\bigl\{B\subseteq\N\mid\exists0\leq r \leq t\dvtx
B^{\{s\}
} \in\underline{\Pi}^{(\infty,0)}_r \\
&&\hspace*{33pt}\mbox{and
contains at least one \emph{white} particle}\bigr\}.
\end{eqnarray*}
Since we started our coloured coalescent in $(\infty,0)$, the
cardinality of $\mathcal B_t$ is at least equal to $A^{\infty}_t$ and,
therefore, we know $\P(\vert\mathcal B_t\vert= \infty)=1$. Since
$\mathcal B_t$ is countable, we can enumerate its elements as $\mathcal
B_t=\bigcup_{n \in\N} \{B^n_t\}$ and use this to define the sets
$\mathcal B^n_t:=\{B^1_t, \ldots, B^n_t\}$, for all $n \in\N$. Since
$\mathcal B_t$ is infinite $\P$-a.s., these $\mathcal B^n_t$ exist for
any $n$, $\P$-a.s.
Now observe that the following inequalities hold even pathwise by construction:
\[
\underline{M}^{(\infty,0)}_t \geq\sum
_{B \in\mathcal B_t} {\mathbf 1}_{\{
B^{\{s\}} \in\underline{\Pi}_t^{(\infty,0)}\}} \geq\sum
_{B \in
\mathcal B^n_t} {\mathbf1}_{\{B^{\{s\}} \in\underline{\Pi}_t^{(\infty
,0)}\}}
\]
and, therefore, the following holds for any $\kappa\in\N$:
\begin{eqnarray*}
\P\bigl( \underline{M}^{(\infty,0)}_t \leq\kappa\bigr) & \leq&\P
\biggl( \sum_{B
\in
B^n_t} {\mathbf1}_{\{B^s \in\underline{\Pi}_t^{(\infty,0)}\}} \leq
\kappa \biggr)
\\
& \stackrel{*} {\leq} &\sum_{i =1}^{\kappa}
\pmatrix{n
\cr
i} \bigl(e^{-ct}\bigr)^i\bigl(1-e^{-ct}
\bigr)^{n-i} \mathop{\longrightarrow}^{n \rightarrow\infty} 0,
\end{eqnarray*}
which in turn implies $ \P( \underline{M}^{(\infty,0)}_t = \infty) =1$.
In * we used that for each of the $n$ blocks in $\mathcal B^n_t$ we
know $\P(B \in\underline{\Pi}^{(\infty,0)}_t) \geq e^{-ct}$ and they
leave the seed-bank independently of each other, which implies that the
sum is dominated by a Binomial random variable with parameters $n$ and
$e^{-ct}$.

Since the probability on the left does not depend on $n$, and the above
holds for any $\kappa\in\N$, we obtain $\P(\underline{M}^{(\infty
,0)}_t = \infty) = 1$ for all $t > 0$. Note that this also implies $\P
(\underline{M}^{(\infty,0)}_t + \underline{N}^{(\infty,0)}_t =
\infty)
= 1$ for all $t > 0$, from which, through the monotonicity of the sum,
we can immediately deduce the stronger statement
\[
\P \bigl(\forall t>0\dvtx\underline{M}^{(\infty,0)}_t +
\underline {N}^{(\infty,0)}_t = \infty \bigr) = 1.
\]
On the other hand, we have seen that $\P(\underline{N}^{(\infty
,0)}_t <
\infty)=1$, for all $t>0$, which again using its monotonicity, yields
$\P(\forall t>0\dvtx\underline{N}^{(\infty,0)}_t < \infty)=1$. Putting
these two results together, we obtain $\P(\forall t>0\dvtx\underline
{M}^{(\infty,0)}_t = \infty) = 1$.
\end{pf*}

\subsection{Bounds on the time to the most recent common ancestor}
\label{ssn:TMRCA}

In view of the previous subsection, it is now quite obvious that the
seed-bank causes a relevant delay in the time to the most recent common
ancestor of finite samples. Throughout this section, we will again use
the notation $(N^{(n,m)}_t, M^{(n,m)}_t)$ to indicate the initial
condition of the block counting process is $(n,m)$.

\begin{defn}
We define the \emph{time to the most recent common ancestor} of a
sample of $n$ plants and $m$ seeds, to be
\[
T_{\mathrm{MRCA}}\bigl[(n,m)\bigr]=\inf\bigl\{t>0\dvtx\bigl(N^{(n,m)}_t,M^{(n,m)}_t
\bigr)=(1,0)\bigr\}. %
\]
\end{defn}

Since coalescence only happens in the plants, $T_{\mathrm{MRCA}}[(n,m)]=\inf\{
t>0\dvtx N^{(n,m)}_t+M^{(n,m)}_t=1\}$. We will mostly be interested in
the case where the sample is drawn from plants only, and write
$T_{\mathrm{MRCA}}[n]:=T_{\mathrm{MRCA}}[(n,0)]$. The main results of this section are
asymptotic logarithmic bounds on the expectation of $T_{\mathrm{MRCA}}[n]$.

\begin{thmm}\label{tmrca}
For all $c,K\in(0,\infty)$, the seed-bank coalescent satisfies
%
\begin{equation}
\E \bigl[T_{\mathrm{MRCA}}[n] \bigr]\asymp{\log\log n}.
\end{equation}
\end{thmm}

Here, the symbol $\asymp$ denotes weak asymptotic equivalence of
sequences, meaning that we have
%
\begin{equation}
\liminf_{n\to\infty}\frac{\E [T_{\mathrm{MRCA}}[n]
]}{\log\log
n}>0
\end{equation}
and
%
\begin{equation}
\limsup_{n\to\infty}\frac{\E [T_{\mathrm{MRCA}}[n]
]}{\log\log
n}<\infty.
\end{equation}

The proof of Theorem~\ref{tmrca} will be given in Propositions \ref
{prop:liminf} and~\ref{prop:upper}.
The intuition behind this result is the following. The time until a
seed gets involved in a coalescence event is much longer than the time
it takes for a plant to be involved in a coalescence, since a seed has
to become a plant first. Thus, the time to the most recent common
ancestor of a sample of $n$ plants is governed by the number of
individuals that become seeds before coalescence, and by the time of
coalescence of a sample of seeds.

Due to the quadratic coalescence rates, it is clear that the time until
the ancestral lines of all sampled plants have either coalesced into
one, or have entered the seed-bank at least once, is finite almost
surely. The number of lines that enter the seed-bank until that time is
a random variable that is asymptotically of order $\log n$, due to
similar considerations as in \eqref{eq:exp_An}. Thus, we need to
control the time to the most recent common ancestor of a sample of
$O(\log n)$ seeds. The linear rate of migration then leads to the
second log.

Turning this reasoning into bounds requires some more work, in
particular for an upper bound. As in the proof of Theorem~\ref
{thmm:notback}, let $X_k, k=1,\ldots,n$ denote independent Bernoulli
random variables with parameters $2c/(k+2c-1)$. Similar to~\eqref
{def:An} define
%
\begin{equation}
\label{def:An_new} A^n:=\sum_{k=2}^n
X_k.
\end{equation}

\begin{lemma}\label{lem:plantseed}
Under our assumptions, for any $\varepsilon>0$,
\[
\lim_{n\to\infty}\P\bigl(A^n\geq(2c+\varepsilon)\log n
\bigr)=0
\]
and
\[
\lim_{n\to\infty}\P\bigl(A^n\leq(2c-\varepsilon)\log n
\bigr)=0.
\]
\end{lemma}

\begin{pf}
As in the proof of Theorem~\ref{thmm:notback} before, we have
\[
\E\bigl[A^n\bigr]=\sum_{k=2}^n
\frac{2c}{k+2c-1}=2c\log n+R'(c,n),
\]
where $R'(c,n)$ converges to a finite value depending on $c$ as $n\to
\infty$, and
\[
\mathbb V \bigl(A^n\bigr)\sim2c \log n\qquad\mbox{as } n\to\infty.
\]
Thus, again by Chebyshev's inequality, for sufficiently large $n$ (and
recalling that~$c$ is our model parameter)
\begin{eqnarray*}
\P\bigl(A^n\geq(2c+\varepsilon)\log n\bigr)&\leq&\P\bigl(A^n-
\E\bigl[A^n\bigr]\geq\varepsilon \log n\bigr)
\\
&\leq&P\bigl(\bigl|A^n-\E\bigl[A^n\bigr]\bigr|\geq\varepsilon\log n
\bigr)
\\
&\leq&\frac{2c}{\varepsilon^2\log n}.
\end{eqnarray*}
This proves the first claim. The second statement follows similarly,
cf. \eqref{eq:chebyshev}.
\end{pf}

Recall the process $(\underline{N}_t, \underline{M}_t)_{t\geq0}$ from
the previous subsection. The coupling of Proposition~\ref{prop:csbc}
leads to the lower bound in Theorem~\ref{tmrca}.

\begin{prop}\label{prop:liminf}
For all $c,K\in(0,\infty)$, the seed-bank coalescent satisfies
%
\begin{equation}
\label{up}\liminf_{n\to\infty}\frac{\E
[T_{\mathrm{MRCA}}[n]
]}{\log
\log n}>0.
\end{equation}
\end{prop}

\begin{pf}
The coupling with $(\underline{N}_t, \underline{M}_t)_{t\geq0}$ yields
\[
T_{\mathrm{MRCA}}[n]\geq\underline{T}_{\mathrm{MRCA}}[n],
\]
where $\underline{T}_{\mathrm{MRCA}}[n]$ denotes the time until $(\underline
{N}_t, \underline{M}_t)$ started at $(n,0)$ has reached a state with
only one block left. By definition, $A^n$ of the previous lemma gives
the number of individuals that at some point become seeds in the
process $(\underline{N}_t, \underline{M}_t)_{t\geq0}$. Thus,
$\underline{T}_{\mathrm{MRCA}}[n]$ is bounded from below by the time it takes
until these $A^n$ seeds migrate to plants (and then disappear). Since
the seeds disappear independently of each other, we can bound
$\underline{T}_{\mathrm{MRCA}}[n]$ stochastically from below
by the extinction time of a pure death process with death rate $cK$
started with $A^n$ individuals. For such a process started at $A^n=l\in
\N$ individuals, the expected extinction time as $l\to\infty$ is of
order $\log l$. Thus, we have for $\varepsilon>0$ that there exists $C>0$
such that
\begin{eqnarray*}
\E \bigl[T_{\mathrm{MRCA}}[n] \bigr]&\geq& \E \bigl[T_{\mathrm{MRCA}}[n]
\mathbf1_{\{A^n\geq(2c-\varepsilon)\log n\}} \bigr]
\\
&\geq& C \log\log n \P \bigl(A^n\geq(2c-\varepsilon)\log n \bigr),
\end{eqnarray*}
and the claim follows from the fact that by Lemma~\ref{lem:plantseed},
$A^n\geq(c-\varepsilon)\log n$ almost surely as $n\to\infty$.
\end{pf}

To prove the corresponding upper bound, we couple $(N_t, M_t)$ to a
functional of another type of coloured process.

\begin{defn}
Let $(\overline{N}_t,\overline{M}_t)_{t\geq0}$ be the continuous-time
Markov process with state space $E\subset\mathbb{N}\times\mathbb{N}$,
characterised by the transition rates:
\[
(n,m)\mapsto\cases{ (n-1, m+1), &\quad$\mbox{at rate } cn$,\vspace *{2pt}
\cr
(n+1,m-1),&\quad $\mbox{at rate } cKm$,\vspace*{2pt}
\cr
(n-1,m),&\quad $\mbox{at rate }
\pmatrix{n
\cr
2} \cdot\mathbf{1}_{\{n\geq
\sqrt{n+m}\}}$. }
\]
\end{defn}

This means that $(\overline{N}_t,\overline{M}_t)_{t\geq0}$ has the
same transitions as $(N_t, M_t)$, but coalescence is suppressed if
there are too few plants relative to the number of seeds. The effect of
this choice of rates is that for $(\overline{N}_t,\overline
{M}_t)_{t\geq0}$, if $n\gtrsim\sqrt{m}$, then coalescence happens at
a rate which is of the same order as the rate of migration from seed to plant.

\begin{lemma}
The processes $(\overline{N}_t,\overline{M}_t)_{t\geq0}$ and $(N_t,
M_t)_{t\geq0}$ can be coupled such that
\[
\P \bigl( \forall t \geq0\dvtx N^{(n,m)}_t \leq
\overline{N}^{(n,m)}_t \mbox{and } M^{(n,m)}_t
\leq\overline{M}^{(n,m)}_t \bigr) =1.
\]
\end{lemma}

\begin{pf}
We construct both processes from the same system of blocks. Start with
$n+m$ blocks labelled from $\{1,\ldots,n+m\}$, and with $n$ of them
carrying an $s$-flag, the others a $p$-flag. Let $S^i, P^i, i=1,\ldots
,n+m$ and $V^{i,j}, i,j=1,\ldots,n+m, i< j$ be independent Poisson
processes, $S^i$ with parameter $cK$, $P^i$ with parameter $c$, and
$V^{i,j}$ with parameter 1. Moreover, let each block carry a colour
flag, blue or white. At the beginning, all blocks are supposed to be
blue. The blocks evolve as follows: At an arrival of $S^i$, if block
$i$ carries an $s$-flag, this flag is changed to $p$ irrespective of
the colour and the state of any other block. Similarly, at an arrival
of $P^i$, if block $i$ carries a $p$-flag, this is changed to an
$s$-flag. At an arrival of $V^{i,j}$, and if blocks $i$ and $j$ both
carry a $p$-flag, one observes the whole system, and proceeds as follows:
\begin{longlist}[(ii)]
\item[(i)] If the total number of $p$-flags in the system is greater or
equal to the square root of the total number of blocks, then blocks $i$
and $j$ coalesce, which we encode by saying that the block with the
higher label ($i$ or $j$) is discarded. If the coalescing blocks have
the same colour, this colour is kept. Note that here the \emph{blocks}
carry the colour, unlike in the coloured process of the previous
sections, where the individuals were coloured. If the coalescing blocks
have different colours, then the colour after the coalescence is blue.
\item[(ii)] If the condition on the number of flags in (i) is not
satisfied, then there is no coalescence, but if both blocks were
coloured blue, then the block ($i$ or $j$) with the higher label is
coloured white (this can be seen as a ``hidden coalescence'' in the
process where colours are disregarded).
\end{longlist}
It is then clear by observing the rates that $(N_t,M_t)$ is equal in
distribution to the process which counts at any time $t$ the number of
blue blocks with $p$-flags and with $s$-flags, respectively, and
$(\overline{N}_t,\overline{M}_t)
$ is obtained by counting the number of $p$-flags and $s$-flags of any
colour. By construction we obviously have $\overline{N}_t\geq N_t$ and
$\overline{M}_t\geq M_t$ for all $t$.
\end{pf}

Define now
\[
\overline{T}_{\mathrm{MRCA}}[m]:=\inf \bigl\{t\geq0\dvtx\bigl(
\overline{N}_t^{(0,m)}, \overline{M}_t^{(0,m)}
\bigr)=(1,0) \bigr\}.
\]

\begin{lemma}\label{lem:uppercoloured}
There exists a finite constant $C$ independent of $m$ such that
\[
\overline{T}_{\mathrm{MRCA}}[m]\leq C\log m.
\]
\end{lemma}

\begin{pf}
Define for every $k\in1,2,\ldots,m-1$ the hitting times
%
\begin{equation}
H_k:=\inf\{t>0\dvtx\overline{N}_t+\overline{M}_t=k
\}.
\end{equation}
We aim at proving that $\E^{0,m}[H_{m-1}]\leq\frac{C}{\sqrt{m}}$
and $\E
^{0,m}[H_{j-1}-H_{j}]\leq\frac{C}{j-1}$ for $j\leq m-1$, for some
$0<C<\infty$. Here and throughout the proof, $C$ denotes a generic
positive constant (independent of $m$) which may change from instance
to instance. To simplify notation, we will identify $\sqrt{j}$ with
$\lceil\sqrt{j}\rceil$, or equivalently assume that all occurring
square roots are natural numbers. Moreover, we will only provide the
calculations in the case of the standard seed-bank coalescent, that is,
$c=K=1$. The reader is invited to convince herself that the argument
can also be carried out in the general case.

We write
$\overline{\lambda}_t$ for the total jump rate of the process
$(\overline{N}, \overline{M})$ at time $t$, that is,
\[
\overline{\lambda}_t=\pmatrix{\overline{N}_t
\cr
2}1_{\{\overline
{N}_t\geq\sqrt{\overline{N}_t+\overline{M}_t}\}}+\overline {N}_t+\overline{M}_t,
\]
and set
\[
\overline{\alpha}_t:=\frac{{\overline{N}_t\choose 2}1_{\{
\overline
{N}_t\geq\sqrt{\overline{N}_t+\overline{M}_t}\}}}{\overline
{\lambda
}_t},\qquad \overline{
\beta}_t:=\frac{\overline{N}_t}{\overline{\lambda}_t}, \qquad\overline{\gamma}_t:=
\frac{\overline{M}_t}{\overline{\lambda}_t}
\]
for the probabilities that the first jump after time $t$ is a
coalescence, a migration from plant to seed or a migration from seed to
plant, respectively. Even though all these rates are now random, they
are well-defined conditional on the state of the process. The proof
will be carried out in three steps.

\emph{Step \textup{1:} Bound on the time to reach $\sqrt{m}$ plants.} Let
%
\begin{equation}
D_m:=\inf\bigl\{t>0\dvtx\overline{N}_t^{(0,m)}
\geq\sqrt {m}\bigr\}
\end{equation}
denote the first time the number of plants is at least $\sqrt{m}$. Due
to the restriction in the coalescence rate, the process $(\overline
{N}^{(0,m)}_t, \overline{M}^{(0,m)}_t)_{t\geq0}$ has to first reach a
state with at least $\sqrt{m}$ plants before being able to coalesce,
hence $D_m<H_{m-1}$ a.s. Hence, for any $t\geq0$, conditional on
$t\leq D_m$ we have $\overline{\lambda}_t=m$ and $\overline
{N}_t<\sqrt
{m}$. Thus, $\overline{M}_t>m-\sqrt{m}$ a.s. and we note that at each
jump time of $(\overline{N}_t, \overline{M}_t)$ for $t\leq D_m$
\[
\overline{\gamma}_s\geq\frac{m-\sqrt{m}}{m}=1-\frac{1}{\sqrt{m}}\qquad
\mathrm{a.s.}\  \forall s\leq t
\]
and
\[
\overline{\beta}_s\leq\frac{1}{\sqrt{m}}\qquad \mathrm{a.s.}\ \forall s\leq
t.
\]
The expected number of jumps of the process $(\overline{N}_t,
\overline
{M}_t)$ until $D_m$ is therefore bounded from above by the expected
time it takes a discrete time asymmetric simple random walk started at
0 with probability $1-1/\sqrt{m}$ for an upward jump and $1/\sqrt{m}$
for a downward jump to reach level $\sqrt{m}-1$. It is a well-known
fact (see, e.g., \cite{Feller}, Chapter XIV.3) that this expectation is
bounded by $C\sqrt{m}$ for some $C\in(0,\infty)$. Since the time
between each of the jumps of $(\overline{N}_t, \overline{M}_t)$, for
$t<D_m$, is exponential with rate $\overline{\lambda}_t=m$, we get
%
\begin{equation}
\label{eq:expDm} \E^{0,m}[D_m]\leq C\sqrt{m}\cdot
\frac{1}{m}=\frac{C}{\sqrt{m}} .
\end{equation}

\emph{Step \textup{2:} Bound on the time to the first coalescence after reaching
$\sqrt{m}$ plants.} At time $t=D_m$, we have $\overline{\lambda}=
{\sqrt{m}\choose 2}+\sqrt{m}+m-\sqrt{m}$, and thus
\[
\overline{\beta}_t=\frac{\sqrt{m}}{{{\sqrt{m}\choose2}}+m}=\frac
{2\sqrt
{m}}{3m-\sqrt{m}}\leq
\frac{C}{\sqrt{m}} \qquad\mbox{a.s.}
\]
and
\[
\overline{\alpha}_t=\frac{m-\sqrt{m}}{3m-\sqrt{m}}\geq\frac
{1}{3}
\biggl(1-\frac{1}{\sqrt{m}} \biggr) \qquad\mbox{a.s.}
\]
Denote by $J_m$ the time of the first jump after time $D_m$. At $J_m$
there is either a coalescence taking place (thus reaching a state with
$m-1$ individuals and hence in that case $H_{m-1}=J_m$), or a
migration. In order to obtain an upper bound on $H_{m-1}$, as a
``worst-case scenario,'' we can assume that if there is no coalescence
at $J_m$, the process is restarted from state $(0,m)$, and then run
again until the next time that there are at least $\sqrt{m}$ plants
(hence after $J_m$, the time until this happens is again equal in
distribution to $D_m$). If we proceed like this, we have that the
number of times that the process is restarted is stochastically
dominated by a geometric random variable with parameter $\frac
{1}{3}(1-\frac{1}{\sqrt{m}})$, and since
\[
\E^{0,m}[J_m-D_m]=\lambda_{D_m}^{-1}=
\frac{1}{{\sqrt{m}\choose
2}+m}\leq \frac{C}{m},
\]
we can conclude [using \eqref{eq:expDm}] that
%
\begin{eqnarray}
\label{Hn} \E^{0,m}[H_{m-1}]&\leq&\E^{0,m}[J_m]
\frac{3\sqrt{m}}{\sqrt
{m}-1}
\nonumber
\\
&= &\bigl(\E^{0,m}[D_m]+\E^{0,m}[J_m-D_m]
\bigr)\frac{3\sqrt
{m}}{\sqrt
{m}-1}
\\
&\leq&\frac{C}{\sqrt{m}}.\nonumber
\end{eqnarray}

\emph{Step \textup{3:} Bound on the time between two coalescences.} Now we want
to estimate $\E^{0,m}[H_{j-1}-H_j]$ for $j\leq m-1$.
Obviously, at time $H_{j}-$, for $j\leq m-1$, there are at least $\sqrt
{j+1} $ plants, since $\overline{N}_t+\overline{M}_t$ can decrease only
through a coalescence. Therefore (keeping in mind our convention that
Gau\ss-brackets are applied if necessary, and hence $\overline
{N}_{H_j}\geq\sqrt{j+1}-1\geq\sqrt{j}-1$ holds) we obtain
$\overline
{N}_{H_j}\geq\sqrt{j}-1$. Hence, either we have $\overline
{N}_{H_j}\geq\sqrt{j}$ and coalescence is possible in the first jump
after $H_j$, or $\overline{N}_{H_j}=\sqrt{j}-1$, in which case
$\overline{\gamma}_{H_j}\geq\frac{j-\sqrt{j}}{j}=1-\frac{1}{\sqrt
{j}}$, meaning that if coalescence is not allowed at $H_j$, with
probability at least $1-\frac{1}{\sqrt{j}}$ it will be possible after
the first jump after reaching $H_j$. Thus, the probability that
coalescence is allowed either at the first or the second jump after
time $H_j$ is bounded from below by $1-\frac{1}{\sqrt{j}}$.

Assuming that coalescence is possible at the first or second jump after
$H_j$, denote by $L_j$ the time to either the first jump after $H_j$ if
$\overline{N}_{H_j}\geq\sqrt{j}$, or the time of the second jump after
$H_j$ otherwise. Then in the same way as before, we see that $\overline
{\alpha}_{L_j}\geq1-\frac{C}{\sqrt{j}}$. Thus, the probability that
$H_{j-1}$ is reached no later than two jumps after $H_j$ is at least
$ (1-\frac{C}{\sqrt{j}} )^2$.
Otherwise, in the case where there was no coalescence at either the
first or the second jump after $H_j$, we can obtain an upper bound on
$H_{j-1}$ by restarting the process from state $(0,j)$. The probability
that the process is restarted is thus bounded from above by $\frac
{C}{\sqrt{j}}$. We know from equation \eqref{Hn} that if started in
$(0,j)$, there is $\E^{0,j}[H_{j-1}]\leq\frac{C}{\sqrt{j}}$. Noting
that $\overline{\lambda}_{H_j}\geq j$, and we need to make at most two
jumps, we have that $\E^{0,m}[L_j]\leq2/j$. Thus, we conclude
%
\begin{eqnarray}
\E^{0,m}[H_{j-1}-H_{j}]&\leq&
\E^{0,m}[L_j] \biggl(1-\frac{C}{\sqrt
{j}}
\biggr)^2 +\frac{C}{\sqrt{j}}\E^{0,j}[H_{j-1}]
\nonumber
\\
&\leq& \frac{2}{j-1} \biggl(1-\frac{C}{\sqrt{j}} \biggr)+ \biggl(
\frac
{C}{\sqrt
{j}} \biggr)^2
\\
&\leq& \frac{C}{j-1}.\nonumber
\end{eqnarray}
These three bounds allow us to complete the proof, since when starting
$(\overline{N}_t, \overline{M}_t)$ in state $(0,m)$ our calculations
show that
%
\begin{eqnarray}
\E \bigl[\overline{T}_{\mathrm{MRCA}}[m] \bigr]&=&\E^{0,m}[H_1]=
\E[H_m]+\sum_{j=2}^{m-1}
\E[H_{j-1}-H_{j}]
\nonumber
\\[-8pt]
\\[-8pt]
\nonumber
&\leq&\frac{C}{\sqrt{m}}+C\sum_{j=2}^{m-1}
\frac{1}{j-1}\sim C\log m
\end{eqnarray}
as $m\to\infty$.
\end{pf}

This allows us to prove the upper bound corresponding (qualitatively)
to the lower bound in \eqref{up}.

\begin{prop}\label{prop:upper}
For $c,K\in(0,\infty)$, the seed-bank coalescent satisfies
%
\begin{equation}
\label{down}\limsup_{n\to\infty}\frac{\E
[T_{\mathrm{MRCA}}[n]
]}{\log
\log n}<\infty.
\end{equation}
\end{prop}

\begin{pf}
Assume that the initial $n$ individuals in the sample of the process
$(\Pi_t^{(n)})_{t\geq0}$ are labelled $1,\ldots,n$. Let
\[
\mathcal S_r:=\bigl\{k\in[n]\dvtx\exists0\leq t\leq r\dvtx k\mbox{
belongs to an $s$-block at time }t\bigr\}
\]
denote those lines that visit the seed-bank at some time up to $t$. Let
\[
\varrho^n:=\inf\bigl\{r\geq0\dvtx\bigl|{\mathcal S}_r^c\bigr|=1
\bigr\}
\]
be the first time that all those individuals which so far had not
entered the seed-bank have coalesced. Note that $\varrho^n$ is a
stopping time for the process $(\Pi_t^{(n)})_{t\geq0}$, and
$N^{(n,0)}_{\varrho^n}$ and $M^{(n,0)}_{\varrho^n}$ are well-defined as
the number of plant blocks, respectively, seed blocks of $\Pi^{(n)}_{\varrho^n}$.
By a comparison of $\varrho^n$ to the time to the most recent common
ancestor of Kingman's coalescent cf. \cite{W09}, $\E[\varrho^n]\leq2$
for any $n\in\N$, and
thus
%
\begin{eqnarray}
\label{eq:bound-kingman} \E \bigl[T_{\mathrm{MRCA}}\bigl[(n,0)\bigr] \bigr]&\leq& 2+\E
\bigl[T_{\mathrm{MRCA}}\bigl[\bigl(N^{(n,0)}_{\varrho
^n},M^{(n,0)}_{\varrho^n}
\bigr)\bigr] \bigr]
\nonumber
\\[-8pt]
\\[-8pt]
\nonumber
&\leq&2+ \E \bigl[T_{\mathrm{MRCA}}\bigl[\bigl(0,N^{(n,0)}_{\varrho^n}+M^{(n,0)}_{\varrho
^n}
\bigr)\bigr] \bigr],
\end{eqnarray}
where the last inequality follows from the fact that every seed has to
become a plant before coalescing. Recall $A^n$ from \eqref{def:An_new}
and observe that
%
\begin{equation}
N^{(n,0)}_{\varrho^n}+M^{(n,0)}_{\varrho^n}\leq
A^n+1 \qquad\mbox {stochastically}.
\end{equation}
This follows from the fact that for every individual, the rate at which
it is involved in a coalescence is increased by the presence of other
individuals, while the rate of migration is not affected. Thus, by
coupling $(N_t, M_t)_{t\geq0}$ to a system where individuals, once
having jumped to the seed-bank, remain there forever, we see that
$N_{\varrho^n}+M_{\varrho^n}$ is at most $A^n+1$.

By the monotonicity of the coupling with $(\overline{N}_t, \overline
{M}_t)$, we thus see from \eqref{eq:bound-kingman}, for $\varepsilon>0$,
%
\begin{eqnarray}
\E \bigl[T_{\mathrm{MRCA}}[n] \bigr]&\leq&2+\E \bigl[\overline
{T}_{\mathrm{MRCA}}\bigl[A^n+1\bigr] \bigr]
\nonumber
\\
&=& 2+\E \bigl[\overline{T}_{\mathrm{MRCA}}\bigl[A^n+1\bigr]{
\mathbf1}_{\{A^n\leq
(2c+\varepsilon)\log n\}} \bigr]
\\
&&{}+\E \bigl[\overline{T}_{\mathrm{MRCA}}\bigl[A^n+1\bigr]{
\mathbf1}_{\{A^n> (2c+\varepsilon
)\log n\}} \bigr].\nonumber
\end{eqnarray}

From Lemma~\ref{lem:uppercoloured}, we obtain
\[
\E \bigl[\overline{T}_{\mathrm{MRCA}}\bigl[A^n+1\bigr]{
\mathbf1}_{\{A^n\leq(2c+\varepsilon
)\log n\}} \bigr]\leq C\log(2c-\varepsilon)\log n\leq C\log\log n,
\]
and since $A^n\leq n$ in any case, we get
\[
\E \bigl[\overline{T}_{\mathrm{MRCA}}\bigl[A^n+1\bigr]{
\mathbf1}_{\{A^n> (2c+\varepsilon
)\log
n\}} \bigr]\leq C\log n\cdot\P\bigl(A^n>(2c+
\varepsilon)\log n\bigr)\leq C.
\]
This completes the proof.
\end{pf}

\begin{rem}
In the same manner as in the proof of Theorem~\ref{tmrca}, one can show
that for any $a,b\geq0$,
\[
\E \bigl[T_{\mathrm{MRCA}}[an, bn] \bigr]\asymp\log \bigl(\log(an)+bn \bigr).
\]
\end{rem}

\begin{appendix}\label{app}\label{sec:appendix}
\section*{Appendix}


\begin{propp}
\label{prop:forgen}Assume $c=\varepsilon N=\delta M$ and $M\rightarrow
\infty$, $N\rightarrow\infty$. Let $(D_{N,M})_{N,M\in\mathbb{N}}$ be
an array of positive real numbers.
Then the discrete generator of the allele frequency process
$(X^N_{\lfloor D_{N,M}t \rfloor},Y^M_{\lfloor D_{N,M}t \rfloor
})_{t\in
\mathbb{R}^+}$ on time-scale $D_{N,M}$ is given by
\begin{eqnarray*}
\bigl(A^Nf\bigr) (x,y)& =& D_{N,M} \biggl[
\frac{c}{N}(y-x)\frac{\partial
f}{\partial x}(x,y) + \frac{c}{M}(x-y)
\frac{\partial f}{\partial
y}(x,y)
\nonumber
\\[-8pt]
\\[-8pt]
\nonumber
& &\hspace*{44pt}{}+ \frac{1}{N}\frac{1}{2}x(1-x)\frac{\partial^2 f}{\partial x^2}(x,y) +R(N,M)
\biggr],
\end{eqnarray*}
where the remainder term $R(N,M)$ satisfies that there exists a
constant $C_1(c,f)\in(0,\infty)$, independent of $N$ and $M$, such that
\[
\bigl|R(N,M)\bigr|\leq C_1\bigl(N^{-3/2}+M^{-2}+N^{-1}M^{-1}+NM^{-3}
\bigr).
\]
\end{propp}

In particular, in the situation where $M=O(N)$ as $N\to\infty$ and
$D_{N,M}=N$ we immediately obtain Proposition~\ref{prop:convergence}.

\begin{pf*}{Proof of Proposition \ref{prop:forgen}}
We calculate the generator of $(X^N_k, Y^M_k)_{k\geq0}$ depending on
the scaling $(D_{N,M})_{N,M \in\N}$. For $f \in\calC^3([0,1]^2)$ we
use Taylor expansion in 2 dimensions to obtain
\begin{eqnarray*}
\bigl(A^Nf\bigr) (x,y) &=& \frac{1}{D_{N,M}} \biggl[
\frac{\partial f}{\partial
x}(x,y) \E_{x,y} \bigl[X^N_1-x
\bigr] + \frac{\partial f}{\partial
y}(x,y) \E_{x,y} \bigl[Y^M_1-y
\bigr]
\\
&&\hspace*{37pt}{} + \frac{1}{2} \frac{\partial^2 f}{\partial x^2}(x,y) \E _{x,y} \bigl[
\bigl(X^N_1-x \bigr)^2 \bigr]
\\
&&\hspace*{37pt}{} + \frac{1}{2} \frac{\partial^2 f}{\partial y^2}(x,y) \E _{x,y} \bigl[
\bigl(Y^M_1-y \bigr)^2 \bigr]
\\
&&\hspace*{37pt}{} + \frac{\partial^2 f}{\partial x\,\partial y}(x,y)\E_{x,y} \bigl[ \bigl(X^N-x
\bigr) \bigl(Y^M_1-y \bigr) \bigr]
\\
&&\hspace*{37pt}{} + \E_{x,y} \biggl[\mathop{\sum_{\alpha,\beta\in\N_0}}_{\alpha
+\beta=
3}
R^{\alpha, \beta}\bigl(X^N_1,Y^N_1
\bigr) \bigl(X^M_1 - x\bigr)^{\alpha
}
\bigl(Y^M_1-y\bigr)^{\beta} \biggr] \biggr],
\end{eqnarray*}
where the remainder is given by
\[
R^{\alpha,\beta}(\bar x,\bar y):= \frac{\alpha+\beta}{\alpha
!\beta
!} \int_0^1
(1-t)^{\alpha+\beta-1}\frac{\partial^3 f}{\partial
x^{\alpha}\,\partial y^{\beta}} \bigl(x - t(\bar x-x),y - t(\bar y-y)
\bigr)\,dt
\]
for any $\bar x, \bar y \in[0,1]$. In order to prove the convergence,
we thus need to calculate or bound all the moments involved in this
representation.

Given $\P_{x,y}$, the following holds: By Proposition~\ref{prop:transprob},
\begin{eqnarray*}
X^N_1&=&\frac{1}{N}(U+Z),
\\
Y^M_1&=&\frac{1}{M}(yM-Z+V),
\end{eqnarray*}
in distribution where $U$, $V$ and $Z$ are independent random variables
such that
\begin{eqnarray*}
U&\sim&\Bin(N-c,x),
\\
V&\sim&\Bin(c,x),
\\
Z&\sim&\Hyp(M,c, yM).
\end{eqnarray*}
Thus, we have
\setcounter{equation}{0}
\begin{equation}
\label{eq:exp} \E_{x,y}[U]=Nx-cx,\qquad \E_{x,y}[V]=cx,\qquad
\E_{x,y}[Z] = cy,
\end{equation}
and moreover
\[
\mathbb V_{x,y}(U)=(N-c)x(1-x).
\]
One more observation is that as $0\leq V\leq c$ and $0\leq Z\leq c$, it
follows that $|Z-cX|\leq c$ and $|V-Z|\leq c$, which implies that for
every $\alpha\in\N$
\begin{eqnarray*}
\bigl|\E_{x,y}\bigl[(Z-cX)^\alpha\bigr]\bigr| & \leq& c^\alpha,
\\
\bigl|\E_{x,y}\bigl[(Z-V)^\alpha\bigr]\bigr| & \leq& c^\alpha,
\end{eqnarray*}
and for every $\alpha,\beta\in\N$
%
\begin{equation}
\bigl|\E_{x,y}\bigl[(Z-cX)^\alpha(V-Z)^\beta\bigr]\bigr|\leq
c^{\alpha+\beta}. \label{cbound2}
\end{equation}
We are now prepared to calculate all the mixed moments needed.
\begin{eqnarray*}
\E_{x,y}\bigl[X^N_1 - x\bigr]&=&
\frac{1}{N}\E_{x,y}[U+Z - Nx]
\\
& =&\frac
{1}{N}\E _{x,y}[U - Nx+cx]+\frac{1}{N}
\E_{x,y}[Z-cx]
\\
& = &\frac{c}{N}(y-x).
\end{eqnarray*}
Here, we used \eqref{eq:exp}, in particular $\E_{x,y}[U - Nx+cx]=\E
_{x,y}[U -\E_{x,y}[U]]=0$. Similarly,
\begin{eqnarray*}
\E_{x,y}\bigl[Y^M_1 - y\bigr] & =&
\frac{1}{M} \E_{x,y}[My+V-Z- My]
\\
&=&\frac{1}{M} \E_{x,y}[V-Z]
\\
& =& \frac{c}{M}(x-y).
\end{eqnarray*}
Noting $X_1^N-x=\frac{1}{N}(U - Nx+cx)+\frac{1}{N}(Z-cx)$ leads to
\begin{eqnarray*}
\E_{x,y}\bigl[\bigl(X^N_1 - x
\bigr)^2\bigr]& = &\frac{1}{N^2} \E_{x,y}\bigl[(U -
Nx+cx)^2\bigr]
\\
&&{}+\frac{2}{N^2} \E_{x,y}[U - Nx+cx]\E_{x,y}[Z-cx]
\\
&&{}+\frac{1}{N^2} \E_{x,y}\bigl[(Z-cx)^2\bigr]
\\
&=& \frac{1}{N^2}\mathbb{V}_{x,y}[U]+\frac{1}{N^2} \E
_{x,y}\bigl[(Z-cx)^2\bigr]
\\
&=& \frac{1}{N}x(1-x)-\frac{c}{N^2}x(1-x)+\frac{1}{N^2}
\E_{x,y}\bigl[(Z-cx)^2\bigr],
\end{eqnarray*}
where
\[
\biggl|-\frac{c}{N^2}x(1-x)+\frac{1}{N^2} \E_{x,y}
\bigl[(Z-cx)^2\bigr] \biggr|\leq\frac
{c^2}{N^2}.
\]
Moreover, we have
\[
\bigl| \E_{x,y}\bigl[\bigl(Y^M_1 - y
\bigr)^2\bigr] \bigr| = \biggl|\frac{1}{M^2} \E _{x,y}\bigl[(V-Z)^2\bigr]
\biggr|\leq\frac{c^2}{M^2}.
\]
Using equation \eqref{cbound2}, we get
\begin{eqnarray*}
&&\bigl|\E_{x,y}\bigl[\bigl(X^N_1-x\bigr)
\bigl(Y^M_1-y\bigr)\bigr] \bigr| \\
&&\qquad\leq \biggl|\frac{1}{NM}\E
_{x,y}[U-xN+cx]\E_{x,y}[V-Z]\biggr |
\\
&&\qquad\quad{}+ \biggl|\frac{1}{NM}\E_{x,y}\bigl[(Z-cx) (V-Z)\bigr] \biggr|
\\
&&\qquad\leq\frac{c^2}{NM}.
\end{eqnarray*}
We are thus left with the task of bounding the remainder term in the
Taylor expansion.
Since $f \in\calC^3([0,1]^2)$, we can define
\[
\tilde{C}^f:=\max\biggl\{\frac{\partial^3 f}{\partial x^{\alpha}\,\partial
y^{\beta}}(\bar x,\bar y)\Big\mid
\alpha, \beta\in\N_0, \alpha+ \beta= 3, \bar x, \bar y \in[0,1]\biggr
\},
\]
which yields a uniform estimate for the remainder in the form of
\[
\bigl\vert R^{\alpha,\beta}(\bar x,\bar y) \bigr\vert\leq\frac{1}{\alpha
!\beta
!\bar{C}^f},
\]
which in turn allows us to estimate
\begin{eqnarray*}
&&\biggl| \E_{x,y} \biggl[\mathop{\sum_{\alpha,\beta\in\N_0}}_{\alpha
+\beta
= 3}
 R^{\alpha, \beta}\bigl(X^N_1,Y^N_1
\bigr) \bigl(X^N_1 - x\bigr)^{\alpha
}
\bigl(Y^M_1-y\bigr)^{\beta} \biggr] \biggr|
\\
&&\qquad \leq\frac{1}{\alpha!\beta!\bar{C}^f}\mathop{\sum_{\alpha
,\beta
\in\N_0}}_{\alpha+\beta= 3}
\E_{x,y} \bigl[\bigl|\bigl(X^N_1 - x
\bigr)^{\alpha
}\bigl(Y^M_1-y\bigr)^{\beta} \bigr|
\bigr].
\end{eqnarray*}
Thus the claim follows if we show that the third moments are all of
small enough order in $N$ and $M$. Observe that for $\alpha\in\{
0,1,2\}
$ we have
%
\begin{equation}
\label{eq:triv_bound} \E_{x,y}\bigl[\bigl|(U-Nx+cx)\bigr|^\alpha\bigr]\leq N.
\end{equation}
For $\alpha=0$, this is trivially true, for $\alpha=1$ it is due to the
fact that the binomial random variable $U$ is supported on $0,\ldots
,N-c$ and $Nx-cx$ is its expectation, and for $\alpha=2$ it follows
from the fact that $(U - N x+cx)^{2}=|(U - N x+cx)^{2}|$ and the
formula for the variance of a binomial random variable.
For $\alpha=3$, it follows, for example, from Lemma~3.1 in \cite
{JFS14} that
%
\begin{equation}
\label{eq:thirdmoment} \E_{x,y}\bigl[\bigl|(U-Nx+cx)\bigr|^3\bigr]=O
\bigl(N^{3/2}\bigr).
\end{equation}
Thus, we get for any $0\leq\alpha, \beta\leq3$ such that $\alpha
+\beta
=3$ that
\begin{eqnarray*}
&&\E_{x,y} \bigl[\bigl|(X_1 - x)^\alpha
\bigl(Y_1^M-y\bigr)^\beta\bigr| \bigr]
\\
&&\qquad=\frac{1}{N^\alpha M^\beta}\sum_{i=0}^\alpha
\pmatrix{\alpha
\cr
i} \E _{x,y}\bigl[\bigl|(U-Nx+cx)^i(Z-cx)^{\alpha-i}(V-Z)^\beta\bigr|
\bigr]
\\
&&\qquad\leq\frac{1}{N^\alpha M^\beta}\sum_{i=0}^\alpha
\pmatrix{\alpha
\cr
i} \E_{x,y}\bigl[\bigl|(U-Nx+cx)^i\bigr|\bigr]
\E_{x,y}\bigl[\bigl|(Z-cx)^{\alpha-i}(V-Z)^\beta\bigr|\bigr]
\\
&&\qquad\leq\frac{1}{N^\alpha M^\beta}\sum_{i=0}^\alpha
\pmatrix{\alpha
\cr
i} N(2c)^{\alpha-i+\beta}1_{\{1,2,3\}}(\alpha)+
\frac
{3(2c)^3}{N^{3/2}}1_{\{3\}}(\alpha)
\\
&&\qquad\leq C \biggl(\frac{1}{NM}+\frac{1}{M^2}+\frac{1}{N^{3/2}}
\frac
{N}{M^3} \biggr),
\end{eqnarray*}
from \eqref{cbound2}, \eqref{eq:triv_bound} and \eqref{eq:thirdmoment},
where the constant $C$ depends only on $c$. This completes the proof.
\end{pf*}
\end{appendix}
%

%





\printaddresses
\end{document}